\documentclass[]{article}
\usepackage{arxiv}
\usepackage[utf8]{inputenc}
\usepackage{amsmath}
\usepackage{hyperref}
\usepackage{cleveref}
\usepackage{amsfonts}
\usepackage{float}
\usepackage{amssymb}
\usepackage{graphicx}
\graphicspath{}
\usepackage[english]{babel}
\usepackage{natbib}
\usepackage{xcolor}
\usepackage{geometry}
\usepackage{authblk}
\usepackage{subcaption}
\usepackage[shortlabels]{enumitem}
\pdfoutput=1

 \newcommand{\ek}{\textcolor{black}} 
\newcommand{\shorttitle}{Nonlinear dimensionality reduction then and now:
AIMs for dissipative PDEs in the ML era}
\newcommand{\vect}[1]{\boldsymbol{#1}}

\title{Nonlinear dimensionality reduction \\
then and now: \\
AIMs for dissipative PDEs in the ML era}
\author[1]{Eleni D. Koronaki}
\author[2]{Nikolaos Evangelou}
\author[3]{Cristina P. Martin-Linares}
\author[4,5]{Edriss S. Titi} 
\author[2,*]{Ioannis G. Kevrekidis}

\affil[1]{Faculté des Sciences, de la Technologie et de la Communication, Université de Luxembourg, Maison du Nombre, Avenue de la Fonte 6, L-4364 Esch-sur-Alzette, Luxembourg}
\affil[2]{Department of Chemical and Biomolecular Engineering and Department of Applied Mathematics and Statistics, Whiting School of Engineering, Johns Hopkins University, 3400 North Charles Street,
Baltimore, MD 21218, USA.}
\affil[3]{Department of Mechanical Engineering, Whiting School of Engineering, Johns Hopkins University, 3400 North Charles Street, Baltimore, USA.}
\affil[4]{Department of Mathematics, Texas A$\And$M University, College Station, TX 77843, USA}
\affil[5]{Department of Applied Mathematics and Theoretical Physics, University of Cambridge, Cambridge CB3 0WA, UK}
\affil[*]{Corresponding author: I.G. Kevrekidis, yannisk@jhu.edu}
\begin{document}

\maketitle

\begin{abstract}
    This study presents a collection of purely data-driven workflows for constructing reduced-order models (ROMs) for distributed dynamical systems. The ROMs we focus on, are data-assisted models inspired by, and templated upon, the theory of Approximate Inertial Manifolds (AIMs); the particular motivation is the so-called post-processing Galerkin method of Garcia-Archilla, Novo and Titi. Its applicability can be extended: the need for accurate truncated Galerkin projections and for deriving closed-formed corrections can be circumvented using machine learning tools. When the right \textit{latent variables} are not \textit{a priori} known, we illustrate how autoencoders as well as Diffusion Maps (a manifold learning scheme) can be used to discover good sets of latent variables and test their explainability. The proposed methodology can express the ROMs in terms of (a) theoretical (Fourier coefficients), (b) linear data-driven (POD modes) and/or (c) nonlinear data-driven (Diffusion Maps) coordinates. Both Black-Box and (theoretically-informed and data-corrected) Gray-Box models are described; the necessity for the latter arises when truncated Galerkin projections are so inaccurate as to not be amenable to post-processing. We use the Chafee-Infante reaction-diffusion and the Kuramoto-Sivashinsky dissipative partial differential equations to illustrate and successfully test the overall framework.
\end{abstract}

\section{Introduction}
Separation of time-scales in dynamical systems is crucial toward the development of Reduced Order Models (ROMs). For a certain class of dissipative evolution equations, the long term dynamics are attracted exponentially fast to smooth invariant objects known as inertial manifolds (IMs), facilitating the construction of ROMs on those. The dynamics on the IM can be described by the Inertial Form (a finite ODE system), which accurately captures the long-term behavior of the original infinite-dimensional system \citep{shvartsman1998nonlinear,jolly1990approximate,titi1990approximate,akram2020priori}. The purpose of this paper is to (somewhat systematically) outline (and demonstrate) links between ``traditional" AIM technology and contemporary data-driven reduction tools, giving rise to ``mathematics-assisted" algorithmic ROM workflows. Such connections had initially been experimentally attempted in the 1990s (e.g.  \cite{KrischerKevrekidisErtl1993,shvartsman2000order}); they are currently experiencing a strong revival due to the explosion in machine-learning-assisted modelling \citep{linot2020deep, anirudh2020improved,lee2020model,bar2019learning,benner2015survey}.

IMs have been proven to exist for only a few systems, and even then, they have not been constructed explicitly \citep{jolly1990approximate}. It is, nevertheless, still possible to find approximations of either the global attractor or the IM itself, i.e. Approximate Inertial Manifolds (AIMs), or the dynamics on it, i.e. the Approximate Inertial Form (AIF), and then track the dynamics in this reduced space. The key ansatz is that the attracting dynamics in the complement space to the AIM, are quickly slaved to, and embodied in, the AIF. Along these lines, the Galerkin projection, as well as nonlinear Galerkin projections on approximate inertial manifolds are also popular choices for reduced order modeling \citep{marion1990nonlinear,jolly1991preserving,shen1990long,benner2015survey}. AIM-based ROMs have been proposed for reaction-diffusion systems \citep{foias1988inertial,adrover2002construction}, the Kuramoto-Sivashinsky equation \citep{foias1988inertial, foias1989exponential,jolly1990approximate}, the two-dimensional Navier-Stokes equations \citep{temam1989induced,temam1989inertial,jauberteau1990nonlinear}, and for the three-dimensional Navier-Stokes equations \citep{guermond2008fully}.

In the late 90s the post-processing Galerkin method was proposed \citep{garcia1998postprocessing,garcia1999postprocessing}, initially in the context of dissipative equations. Post-processing Galerkin takes into account the observation that the error between the result of integrating a truncated Galerkin on the one hand, and the projection of the true solution in the finite-dimensional Galerkin space, on the other, is significantly smaller than the error between the truncated Galerkin and the full solution (superconvergence \citep{garcia1999approximate,wahlbin2006superconvergence}). We will return to this below and illustrate it in Sec. \ref{sec:Results} and Fig. \ref{fig:TotalFigure}. Given this observation, one uses the dynamics expressed only in terms of the leading low modes (a truncated version of the equations) to integrate. Once the time integration is finished one can \textit{post-process} the obtained solution by approximating the high modes as a function of the solution in the leading modes. Since, in the post-processing Galerkin framework, the correction is computed only at the end of time integration, this makes it much cheaper to implement computationally than true nonlinear Galerkin \citep{garcia1998postprocessing,garcia1999postprocessing}. Moreover, truncation analysis derivation of the spectral method for dissipative evolution equations, such as the Navies-Stokes equations, gives rise to the \textit{post-processing} Galerkin as the \textit{leading order numerical scheme}, and not the Galerkin scheme itself, as commonly believed \citep{margolin2003postprocessing}.

Model identification assisted by machine-learning has emerged in the 90s and is now experiencing a rebirth as a tool to discover minimal parametrizations of an IM, which can subsequently be used to evolve the dynamical system in a reduced space \citep{lu2017data,chorin2015discrete,zeng2023autoencoders,zeng2022data,de2023data,linot2023turbulence}. Some efforts implemented linear methods like POD \citep{KrischerKevrekidisErtl1993, kang2015nonlinear,shvartsman2000order}, to identify a suitable subspace that contains the majority of the variance of the system, and parametrizes the long term dynamics. More recently, operator inference with quadratic manifolds has been proposed for model reduction \citep{geelen2023operator,zastrow2023data,qian2022reduced,mcquarrie2021data}. Nonlinear dimensionality reduction methods, such as autoencoders \citep{kramer1991nonlinear} or Diffusion Maps (DMAPs) \citep{r21} have also been used to discover latent variables of data that originally live in a high-dimensional space. Learning a dynamical systems in the latent space of an autoencoder (even as a collection of local charts), or in Diffusion Maps space, also provides a systematic approach to ROM construction (e.g. \cite{rico1992discrete,sonday2010manifold,evangelou2023double,linot2022data,lee2020model, bar2019learning}). Needless to say, nonlinear system identification assisted by machine learning remains a very active current research endeavor, encompassing a plethora of directions from symbolic methods e.g. \citep{brunton2016discovering}, to physics-informed  methods, e.g.\citep{raissi2019physics}, to numerics-informed methods \citep{bar2019learning}.

In our view, the ``1980s" IM and AIM efforts towards useful reduced order models of dissipative PDEs can be succinctly summarized as follows:
Given the functional form of the PDE for which we know (or believe) an IM exists, and having an estimate of the dimensionality of said manifold:
\begin{enumerate}[(A)]
    \item start by finding the (leading) eigenmodes, say $k$ of them, of the (dissipative part of the) operator that ``determines" (parametrizes) the IM. In that sense, the components of the solution in the remaining ``higher order" eigenmodes can be expressed as functions of the components in the lower, determining, ones;
    \item guided by separation of time scales ideas, construct the AIM approximating this function, by writing the components of the higher eigenmodes as (approximate) functions of the components of the lower, determining ones. Several implementable such approximations have been proposed and analysed: e.g. the ``steady" manifold, the ``Euler-Galerkin", and the Foias-Manley-Temam (FMT) manifold among others. We already have a practical result: if somebody provides \textit{as observations} the lower mode amplitudes, we can meaningfully \textit{and analytically} improve the full spatiotemporal solution, complementing it with the recovered higher mode components. We will return to this theme when discussing post-processing Galerkin. Let it be noted here that even though the original motivation of AIM was to find an approximation to the IM whenever the latter exists, however, this idea was generalized later and implemented by finding a manifold which approximates the global attractor as a set; observing that global attractor always exists for genuine dissipative dynamical system;
    \item beyond just correcting such observations, these functions can be used to correct approximations of the dynamics through their low-order Galerkin truncation: From an accurate, high-order Galerkin truncation, we keep only the low, ``determining" Galerkin ODEs; instead of omitting the higher order terms as negligible,  we now  {\em substitute} the AIM function in the low terms. We now have the ``steady", or ``Euler-Galerkin" or FMT inertial forms.
\end{enumerate}

This original program is complemented by the ``post-processing Galerkin" protocol: Here we actually keep the low order Galerkin truncation, ignoring the contribution of the higher order, slaved modes to it, expecting/believing that, in its low-dimensional space, these few ODEs are accurate enough to approximate the projection of the exact solution on the Galerkin space. The authors of \citep{garcia1998postprocessing,garcia1999postprocessing,garcia1999approximate} took into account the observation that the total error of the solutions predicted by the truncated low-order Galerkin is appreciably larger than the error after adding to them (in a sense, ``reinjecting") the AIM-approximated higher order solution components. This reinjection
is performed {\em after} the truncated low-order Galerkin equations have been integrated until each time instance of interest (we remind the reader that this is revisited in Sec. \ref{sec:Results} and Figure \ref{fig:TotalFigure}).  

They named the approach ``post-processing Galerkin" since it takes place after the truncated low-order Galerkin has been obtained and integrated: it is these concrete available {\em solutions} of the model that are being improved - not the model itself.

Explicit AIMs have been obtained in the context of spectral Galerkin approximation by writing approximations of the evolution equation of the high modes in terms of the low modes, a closure relation. In the context of spectral Galerkin approximation based on Fourier modes or eigenfunctions of the Stokes operator, one can naturally decouple the phase space into low Fourier (eigenfunctions) modes and their complement high Fourier (eigenfunctions) modes. Therefore, the above-described strategy of obtaining AIM is possible to be executed explicitly, and leads to an analytical closure. For the examples we present in this work, the spectral Galerkin approximation could indeed provide a desirable closure.

However, we would like to note that in the context of the Finite Element Galerkin method, the above decomposition to coarse spatial scales and their complement is not a straightforward task. Therefore the above strategy can not be followed to obtain an explicit (paper and pencil) closure form, that expresses the fine spatial scales of the solution in terms of the coarse finite elements spatial scales. For this case, a more general framework for implementing the Post-processing Galerkin can be used \citep{garcia1999postprocessing}. In this more general case, an explicit form of an AIM in order to implement Post-processing Galerkin is not required. We briefly present this more general scheme in Sec \ref{sec:PPE_FE}.

Today, beyond symbolic model (AIM) or solution (post-processing Galerkin) improvement, data driven techniques allow us (given accurate simulation data or observations) to:
\begin{enumerate}[(a)]
    \item Estimate the AIM dimensionality in a data-driven way (either through autoencoders or through manifold learning).
    \item Learn good reduced AIFs (the ``correct'', nonlinear Galerkin, right hand-side of the reduced, low order, components of the PDE) in a data-driven
way.
    \item Learn the AIM functions (high order mode components as a function of low order model components) in a data driven way.
    \item Given the learned AIM in (c), correct the solutions of a low-order Galerkin truncation (a ``data driven" post-processing Galerkin). Beyond the steps (b-d) above, that more or less correspond to the traditional (A-C) analytical steps, there are now a couple of very useful data-driven ``twists".
    \item Circumvent the assumption of accuracy of the low-order (linear) Galerkin truncation; the low order AIF is learned from observations of the low-order components of accurate full PDE dynamics; and now the ``post-processing" that follows can be done (1) with the same ``old" analytical AIMs, or, interestingly (2) with data-driven learned AIMs from the same accurate full PDE dynamics.
    \item Gray-Box (in some sense ``physics-assisted") learning: instead of a fully black-box learning of the AIF using PDE observations, we now learn {\em the correction} of the not-so-quantitative low-order linear Galerkin truncation. This correction can be learned as an additive (residual) term, or even as a {\em functional} correction - hoping for easier training, since what is learned is a perturbation of the identity \citep{martin2023physics}.
    \item (\textit{This is not so much a step in our list, as a branching towards new capabilities}). Up to now, everything {\em but the eigenfunctions} parametrizing the manifold was data-driven; the eigenfunctions themselves were still analytical. If we allow ourselves to find the parametrization of the manifold in a data-driven way, two individually significant new options arise:
    \begin{enumerate}
        \item Use linear data-driven eigenfunctions: the leading Principal Components (PODs) of the full accurate PDE simulations.  Now the low-order PODs parametrize the manifold, and the higher order POD components embody the AIM. POD-Galerkin takes the place of traditional Galerkin.
        \item Use a \textit{nonlinear}, data-driven AIM parametrization: One can here either (g2a) use the latent variables of an autoencoder to parametrize the AIM, learn the corresponding accurate AIF, {\em and} post-process it to more accurate spatiotemporal PDE solution reconstruction; or (g2b) use the leading POD components to parametrize the AIM, learn the accurate corresponding AIF, {\em and } post-process it for a more accurate spatiotemporal PDE solution reconstruction. At the risk of making this list ridiculously long, we also add -and illustrate below- the possibility (g2c) of using spectral (Diffusion Map) data mining to parametrize the AIM along with the associated \textit{Geometric Harmonics} for the post-processing.
    \end{enumerate}

\end{enumerate}

\begin{figure} [ht!]
    \centering
    \includegraphics[width=1\textwidth]{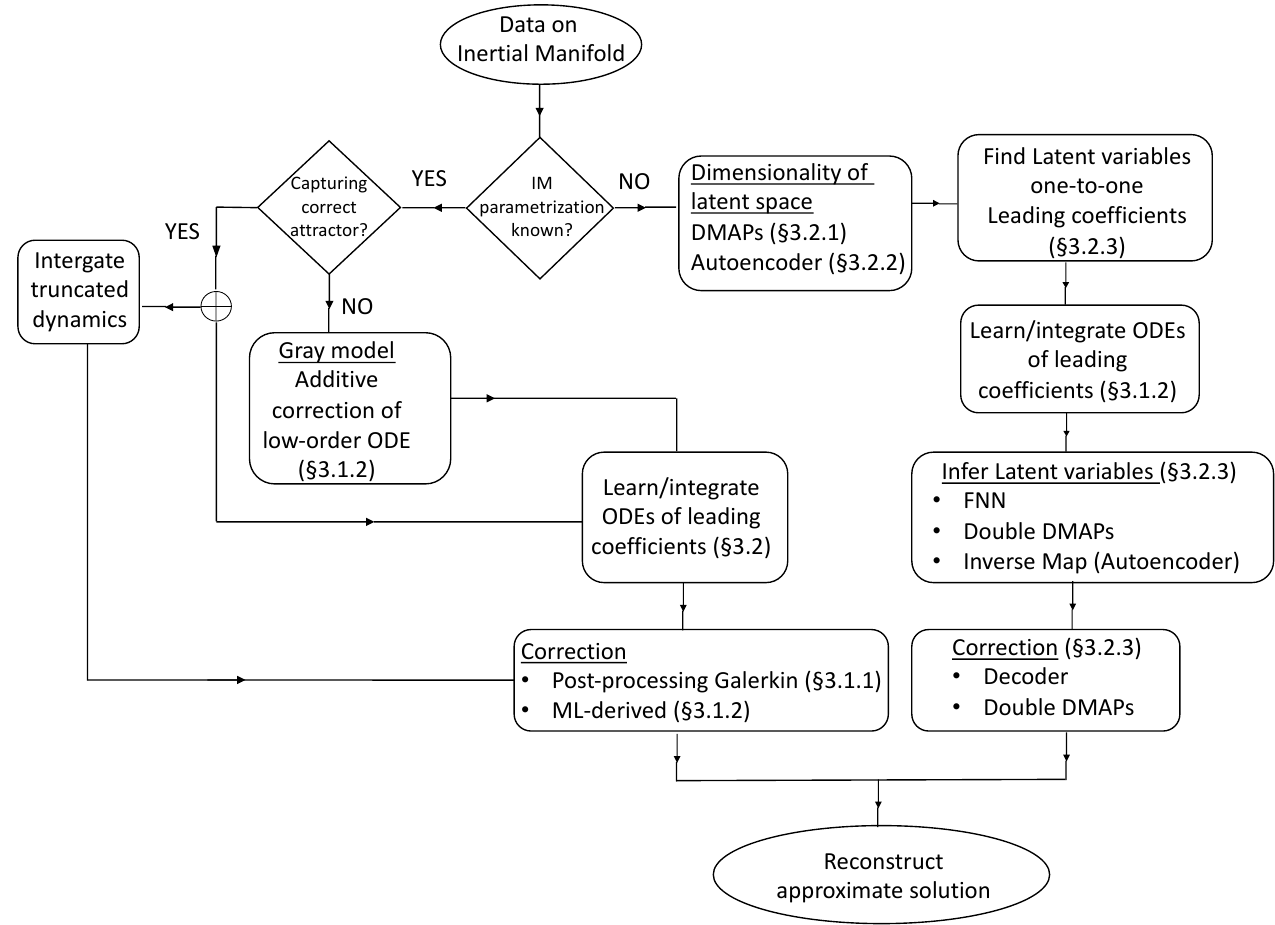}
    \caption{Flowchart of the proposed workflow}
    \label{fig:Overview}
\end{figure}

A schematic overview of the different options proposed in each case, are presented in Fig.\ref{fig:Overview}, with references to the subsequent sections where they are discussed in detail.

The remainder of the paper is organized as follows: After listing the illustrative examples used in this study (Sec \ref{sec:Iluexamp}), we proceed with describing the methodology (Sec. \ref{sec:method}.) We start with briefly reviewing the ``traditional" approximations of IMs and IFs (and AIMs and AIFs) (Sec. \ref{subsec:LearnHarm}). We then discuss neural network-based alternatives to approximating IMs and IFs (Sec. \ref{sec:NNAIMandAIF}), followed by nonlinear manifold learning methods for determining the dimensionality and parametrization of the latent space (Sec. \ref{sec:dmaps} and \ref{sec:autoencoders}). 
After presenting our results we conclude by pointing out that the technology can be easily ``transferred" to POD parametrizations of the IM (Sec. \ref{sec:AIMtransf}).

\section{Illustrative examples: The Chafee-Infante and the Kuramoto-Sivashinsky equations}
\label{sec:Iluexamp}
Our first example is the reaction-diffusion Chafee-Infante partial differential equation (PDE), for which the dimensionality of the Inertial Manifold (IM), for the parameter range of interest, is known; it reads:
\begin{equation}
\label{eq:Chafee_Infante}
    u_t = u - u^3 + \nu u_{xx}.
\end{equation}

\noindent{The} parameter $\nu$ was chosen as $\nu = 0.16$ and Dirichlet boundary conditions, $u(0,t) = u(\pi,t)=0$, were used. The Chafee-Infante PDE, for $\nu=0.16$, has been shown to have a two-dimensional inertial manifold \citep{sonday2011systematic,gear2011slow,jolly1989explicit,evangelou2022double}. To simulate the dynamics on/near this two-dimensional manifold, the Galerkin projection
\begin{equation}
    u(x,t) \approx \sum_{k=1}^{3} \alpha_{k}(t)sin(kx)
\end{equation}
 was used \citep{,gear2011slow,jolly1989explicit,evangelou2022double}. The first two leading sine coefficients {$\alpha_1(t), \alpha_2(t)$} are sufficient to parameterize this two-dimensional manifold, and Galerkin equations based only on the first two modes provide a qualitatively correct approximation of the dynamics in these two modes. We will consider the solution of the Chafee-Infante equation with three modes as the ground truth (cf Fig.\ref{fig:ChaInf2D3D}a). For the post-processing process, the truncated equations with the first two sine coefficients $\tilde{\vect{\alpha}} = \{\alpha_1, \alpha_2 \}$ are the ones used for integration up to time $t= T$. Then, their solution is post-processed to recover $\hat{\alpha} = \alpha_3$ and reconstruct the full solution $u(x,T)$. For this first example, the truncated dynamics governed by the first two sine coefficients are (considered to be) qualitatively, but not quantitavely accurate; the post-processing step aims to correct the obtained solution from these truncated dynamics.

To demonstrate the potential of the proposed methodology in a case with more complex dynamics, we select the Kuramoto-Sivashinsky (KS) PDE,
\begin{equation}
   \begin{split}
    & u_t=-\nu(u u_x+ u_{xx})-4u_{xxxx}; \;\; \textrm{for} \;\; x\in[0,2\pi].\\
  \end{split}
   \label{eq:KS}
\end{equation}

\noindent{The KS} (Equation \eqref{eq:KS}) is a prototypical equation with dynamics that include chaos, derived in the context of a diverse range of physical systems such as, but not limited to, thin film flow on inclined planes and instabilities in a laminar flame front \citep{kuramoto1976persistent,sivashinsky1977nonlinear,alekseenko1985wave,chang1986nonlinear,chang1986traveling,jolly1990approximate,Kevrekidis1990}. 
The parameter $\nu$ in our case is set to $\nu=33$ and periodic boundary conditions are used $u(0,t)= u(2\pi, t)$.
In this example, Fourier series expansion with 8 terms is used to approximate the ground truth  $u(x,t)$:
\begin{equation}
u(x,t)\approx\sum_{k=1}^{8} \alpha_k(t)sin(kx)+ \beta_k(t) cos(kx); \;  x \in [0,2\pi],
\end{equation}
 \noindent 
 which results in 8 ODEs for the sine coefficients ($\{\alpha_k\}_{k=1}^8$ ) and 8 for the cosine coefficients ($\{\beta_k\}_{k=1}^8$).

 Restriction to the space of odd functions leads to retaining only the sine terms, resulting in a system of 8 ODEs for the sine coefficients which is considered, in this work, as the exact solution of the KS. We use the truncation to the leading three sine coefficients $\tilde{\vect{\alpha}} = \{\alpha_1,\alpha_2,\alpha_3\}$ to study the  dynamics for $\nu=33$; however even though it has been shown that a 3D manifold exists, the truncated equations based on the leading coefficients do not provide an accurate approximation of the dynamics of these coefficients. In this case the \textit{traditional} post-processing Galerkin methodology, does not apply (we do not have a good base solution to correct). We circumvent this issue by constructing Gray-Box models, as we show below in Sec. \ref{subsub:trunc}.

\section{Methodology}
\label{sec:method}
\subsection{Approximating the IM and the IF (known latent space)}
\label{subsec:LearnHarm}
\subsubsection{Euler-Galerkin}
\label{sec:Euler_Galerkin}


As a preamble to  \textit{traditional} post-processing Galerkin, here we discuss {\em nonlinear} Galerkin schemes, in particular the ``Euler-Galerkin'' algorithm, that provides a closed-form approximation of inertial manifolds \citep{foias1988computation}. Consider the evolution equation

\begin{equation}
\label{eq:general_evolution}
    \frac{du}{dt} + Au + F(u) = 0, \hspace{0.5cm} u\in H
\end{equation}

\noindent{where} $H$ is an appropriate Hilbert space, $A$ is a self-adjoint positive-definite linear operator with compact inverse, and let $F$ be a nonlinear operator such that equation \eqref{eq:general_evolution} is globally well-posed in time for all initial data in $H$.
By denoting a projection onto the span of the first $n$ eigenvectors of $A$ by $P$ and $Q = I -P$ we can split Equation \eqref{eq:general_evolution} into
\begin{equation}
       \frac{dp}{dt} + Ap + PF(p+q) = 0
\end{equation}
\begin{equation}
\label{eq:dq}
       \frac{dq}{dt} + Aq + QF(p+q) = 0
\end{equation}

where $p=Pu$, $q=Qu$ and $q+p = u$. Assuming that the long-term dynamics of Equation \eqref{eq:general_evolution} live in a $n-$dimensional inertial manifold described as the graph of a function $\Phi: PH \to QH$ we can write the projection of the inertial manifold onto $PH$ as 

\begin{equation}
    \frac{dp}{dt} + Ap + PF(p + \Phi(p)) =0.
\end{equation}

An approximation of $\Phi$ is achieved through a Galerkin truncation of $m$ modes in Equation \eqref{eq:dq} , where $m>n$. The projection to the space of the higher modes $n+1,\dots,m$ defines $Q_m$. Since the higher modes are attracted exponentially fast to the IM and become functions of the lower modes, we perform an implicit Euler step to approximate the solution $\hat{q}$ with a step size $\tau$. By assuming an initial condition $q_0=0$ we get

\begin{equation}
\label{eq:euler_step1}
    \hat{q} + \tau A\hat{q} + \tau Q_mF(p + \hat{q}) = 0.
\end{equation}

Instead of completely solving equation \eqref{eq:euler_step1} we perform a single fixed-point iteration using an initial $\hat{q}=0$ and holding the lower modes $1, \dots, n$ (the components of $p$) constant. This gives the approximation:

\begin{equation}
    \hat{\Phi}_m (p) = -\tau(I + \tau A)^{-1} Q_mF(p)
\end{equation}

\noindent{an algebraic} expression that estimates the higher modes $\{n+1, \dots, m\}$ as a function of the lower $n$ modes and thus an approximation of the IM itself.

Substituting $\hat{\Phi}_m (p)$ for the $m-n$ higher modes gives

\begin{equation}
    \label{eq:Substituting}
        \frac{dp}{dt} + Ap + PF(p + \hat{\Phi}_m (p)) =0
        \end{equation}
  and more precisely an Euler-Galerkin approximation consisting of $n$ differential equations 
\begin{equation}
        \label{eq:Substituting2}
        \frac{dp}{dt} + Ap + PF(p -\tau(I + \tau A)^{-1} Q_mF(p)) =0.
\end{equation}
  
In this work, the (nonliner) Euler-Galerkin algorithm was applied to the Chafee-Infante partial differential equation, as detailed in Sec. \ref{sec:Euler_Garlerkin_Chafee_Infante}.
\subsubsection{Neural network derived AIM and AIF}
\label{sec:NNAIMandAIF}

\begin{figure}[ht!]
    \centering
    \includegraphics[width=0.3\textwidth]{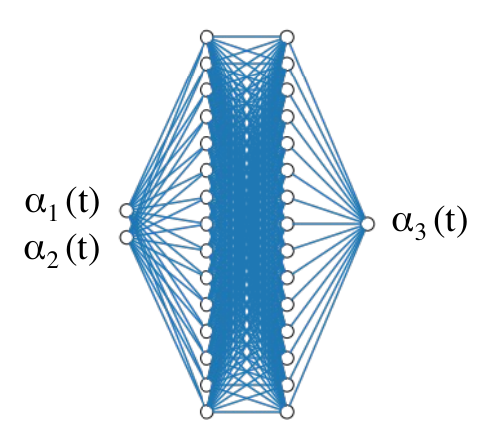}
    \caption{Illustrative example of a feed-forward neural network for prediction of higher coefficients. In this example the lower harmonics, $\alpha_1(t)$ and $\alpha_2(t)$ are used as inputs to the network that predicts $\alpha_3(t)$}
    \label{fig:NNhigher}
\end{figure}

The higher sine modes' coefficients $\hat{\vect{\alpha}}$,  which are necessary for accurate reconstruction of the solution in physical space, can be obtained in a data-driven manner. Specifically, here we use deep neural networks, schematically shown in Fig.\ref{fig:NNhigher}, to learn the coefficients $\hat{\vect{\alpha}}$, given the values of leading (lower) sine modes' coefficients at a specific point in time, $t=T$, $\tilde{\vect{\alpha}}(T)$. 
\begin{equation*}
\hat{\vect{\alpha}}(t)=f_{NN}(\tilde{\vect{\alpha}}(t)), 
\end{equation*}
where $\tilde{\vect{\alpha}}$ stands for leading sine coefficients (low modes) and $\hat{\vect{\alpha}}$ stands for the higher sine modes coefficients. The leading coefficients $\hat{\vect{\alpha}}(T)$ have been obtained as a result of the time integration of the truncated dynamics.

Alternatively, when the result of time-integration of the two truncated lower sine coefficients equations, is inaccurate, we can correct it by learning a data-driven truncated
 ODE in the lower sine coefficients, with general form:
\begin{equation}
    \frac{d\vect{\alpha}}{dt}=\vect{f}(\vect{\alpha})
    \label{eq:LearnedODE}
\end{equation}
\noindent
where $\vect{\alpha} \in \mathbb{R}^m$, here $m=2$, are the variables in which we observe the evolution of the dynamics. Observe that since: m=2 here the Poincaré-Bendixon theorem applies. Hence the dynamics of the low modes is either goes to a limit-cycle or to a steady state. This data-driven AIF was first explicitly described and implemented in \citep{shvartsman2000order} (see also in \citep{KrischerKevrekidisErtl1993}).

\begin{figure}[ht!]
    \centering
    \includegraphics[width=0.25\textwidth]{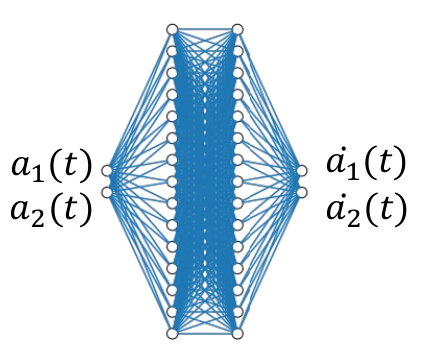}
    \caption{An example of a feed-forward neural network architecture for the approximation of the right-hand-side of an evolution ODE.} 
    \label{fig:LearnODE}
\end{figure}

The function $\vect{f}$ is approximated by a fully connected neural network, schematically represented in Fig.\ref{fig:LearnODE}. The goal is to predict the time derivatives of the lower sine coefficients from the their values. Once this is done, the right-hand-side of the ODEs in Eq.~\ref{eq:LearnedODE} can be used in conjunction with any method of integration in time, such as the Runge-Kutta, to accurately approximate $\tilde{\alpha}(T)$ and then proceed as above to post-process $\tilde{\alpha}(T)$.

For parameter values of the KS equation for which the long-term truncated dynamics may not be accurate, an appealing alternative to the Black-Box approach discussed above arises. 

One can remedy the situation by \textit{first correcting} the reduced dynamics, before deriving the missing terms for reconstruction. This can be achieved by constructing a ``Gray-Box" data-driven dynamic model. This Gray-Box model describes the evolution of a reduced system, by adding to the truncated dynamics a \textit{learned} correction term, which can be thought of as a closure. This correction is approximated by a neural network that takes as inputs the lower order sine coefficients and delivers the difference between their true time-derivatives and the truncated Galerkin time-derivatives:
\begin{equation*}
\frac{d{\tilde{\vect{\alpha}}}^{\text{p}}}{dt}-\frac{d{\tilde{\vect{\alpha}}}^{\text{{t}}}}{dt} =g_{NN}(\tilde{\vect{\alpha}}(t)) 
\end{equation*}
\noindent{where} $\frac{d{\tilde{\vect{\alpha}}}^{\text{p}}}{dt}$ is the true vector field projected in the leading sine coefficients $\tilde{\vect{\alpha}}$ and $\frac{d{\tilde{\vect{\alpha}}}^{\text{t}}}{dt}$ is the vector field of the corresponding truncated Galerkin projection.  

Here, $g_{NN}$ is approximated using a neural network implemented in tensorflow \citep{tensorflow2015-whitepaper} with 6 hidden layers, 95 neurons each, and a $tanh$ activation function. The loss function used is the mean squared error (MSE), and the Adam optimizer is employed.


Finally, it is worth noting that the proposed workflow works equally well, when considering evolution equation of the leading POD mode coefficients, as parametrizing the IM. An illustrative example, based on the Chafee-Infante POD-based equations can be found in \ref{sec:Euler_Garlerkin_Chafee_Infante}.

\subsection{Learning the dimensionality of the latent space}
\label{sec:LearnLatent}
\begin{figure}[ht]
    \centering
    \begin{subfigure}[b]{0.45\textwidth}
         \centering
         \includegraphics[trim={0.0cm 0.5cm 0.0cm 0.0cm},clip,width=\textwidth]{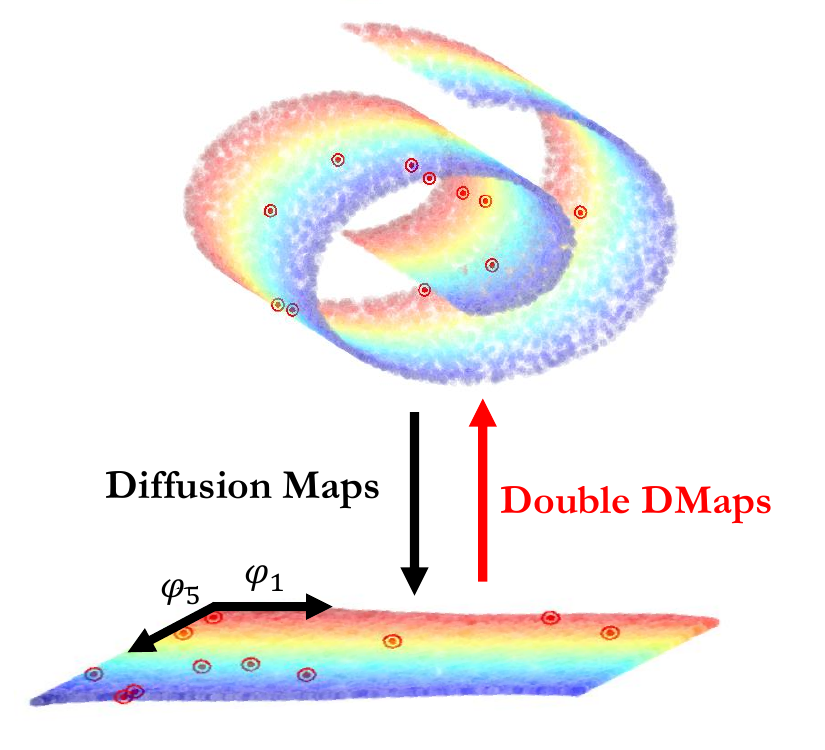}
         \caption{}
    \end{subfigure}
    \hfill
    \begin{subfigure}[b]{0.45\textwidth}
         \centering
         \includegraphics[width=\textwidth]{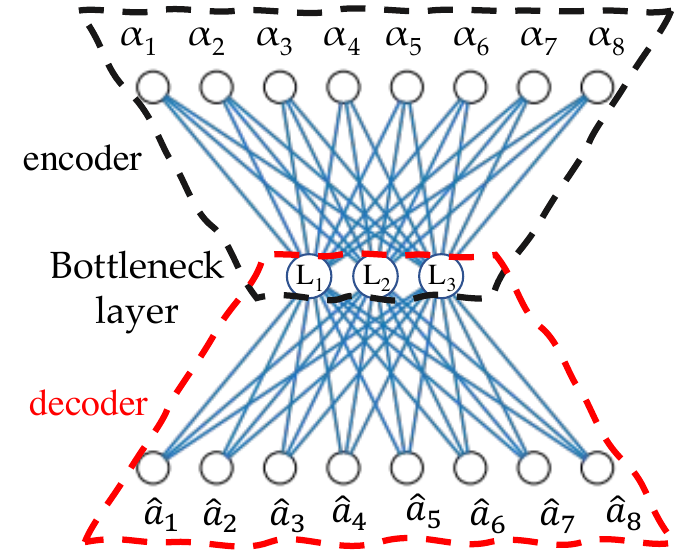}
         \caption{}
    \end{subfigure}
    \caption{Learning a low dimensional embedding of data: (a) Manifold learning with Diffusion Maps and inverse transformation with Double Diffusion Maps \citep{evangelou2022double}; (b) Representative autoencoder structure, including encoder/decoder and the bottleneck layer.}
\label{fig:ManifoldLearning}
\end{figure}

In most cases, a minimal parametrization of the IM of a dynamical system is not known \textit{a priori}. It is possible to discover it, using different data mining approaches, such as Diffusion Maps and autoencoders. Both methods are discussed in the following paragraphs and summarized in Fig.\ref{fig:ManifoldLearning}

\subsubsection{Diffusion Maps}
\label{sec:dmaps}
Diffusion maps ~\citep{r19,r20,r21} is a manifold learning framework that can (based upon diffusion processes) facilitate discovering low-dimensional intrinsic geometric descriptions of data sets, even when the data is high-dimensional, nonlinear and/or corrupted by (relatively small) noise. It is used here to discover the dimensionality of the IM and provide a data-driven parametrization of it.

The parametrization of the manifold is obtained through
a few eigenvectors, $\phi_i$, of a scaled affinity matrix, which contains the Euclidean distances between all the pairs of available data points. A detailed description of the Diffusion Maps algorithm is provided in the Sec.  \ref{sec:Diffusion_Maps_algorithm} of the Appendix and for the Double Diffusion Maps in Sec. \ref{sec:double_dmaps}.

\subsubsection{Autoencoders}
\label{sec:autoencoders}

Autoencoders \citep{kramer1991nonlinear} are neural networks that are trained (a) to encode high-dimensional data into a low-dimensional representation (b) to reconstruct the original high-dimensional from this lower-dimensional representation (cf.Fig.\ref{fig:ManifoldLearning}b). In this context, the input layer is the same as the output, and the low-dimensional encoding is parametrized by the weights of the bottleneck layer. The loss function 
\begin{equation}
    L_{\theta} = || \vect{\alpha}^{(k)}-\bar{\vect{\alpha}}^{(k)} ||^2
\end{equation}
\noindent{is} commonly used to train an autoencoder where $\vect{\alpha}^{(k)}$ represent a data point in the ambient space and $\bar{\vect{\alpha}}^{(k)}$ the reconstructed data point $k$ from the autoencoder.

In this work, we use autoencoders for an additional second use case, which relies on the observation that the discovered autoencoder latent coordinates are one-to-one with the leading sine coefficients  $\tilde{\vect{\alpha}}$, as discussed in detail in the following paragraph.

\subsubsection{Theoretical and data-driven latent variables: transformations and AIMs}
\label{sec:AIMtransf}
The \textit{local} one-to-one relation between the autoencoder's latent variables ($\vect{L}$) and the leading sine coefficients ($\vect{\tilde{\alpha}}$) is tested by computing the Inverse Function Theorem across the training data. The Inverse Function Theorem guarantees local invertibility in a neighborhood of any point $\vect{L}_i \in \vect{L}$ if the determinant of the Jacobian (det($\mathbf{J}_f(\mathbf{L}))$ is bounded away from zero. We provide a more detailed description of the Inverse Function Theorem in Sec. \ref{sec:IFT} of the Appendix. The Jacobian computation in our case is performed by using automatic differentiation with tensorflow.

\begin{figure}[ht!]
    \centering
    \includegraphics[width=0.8\textwidth]{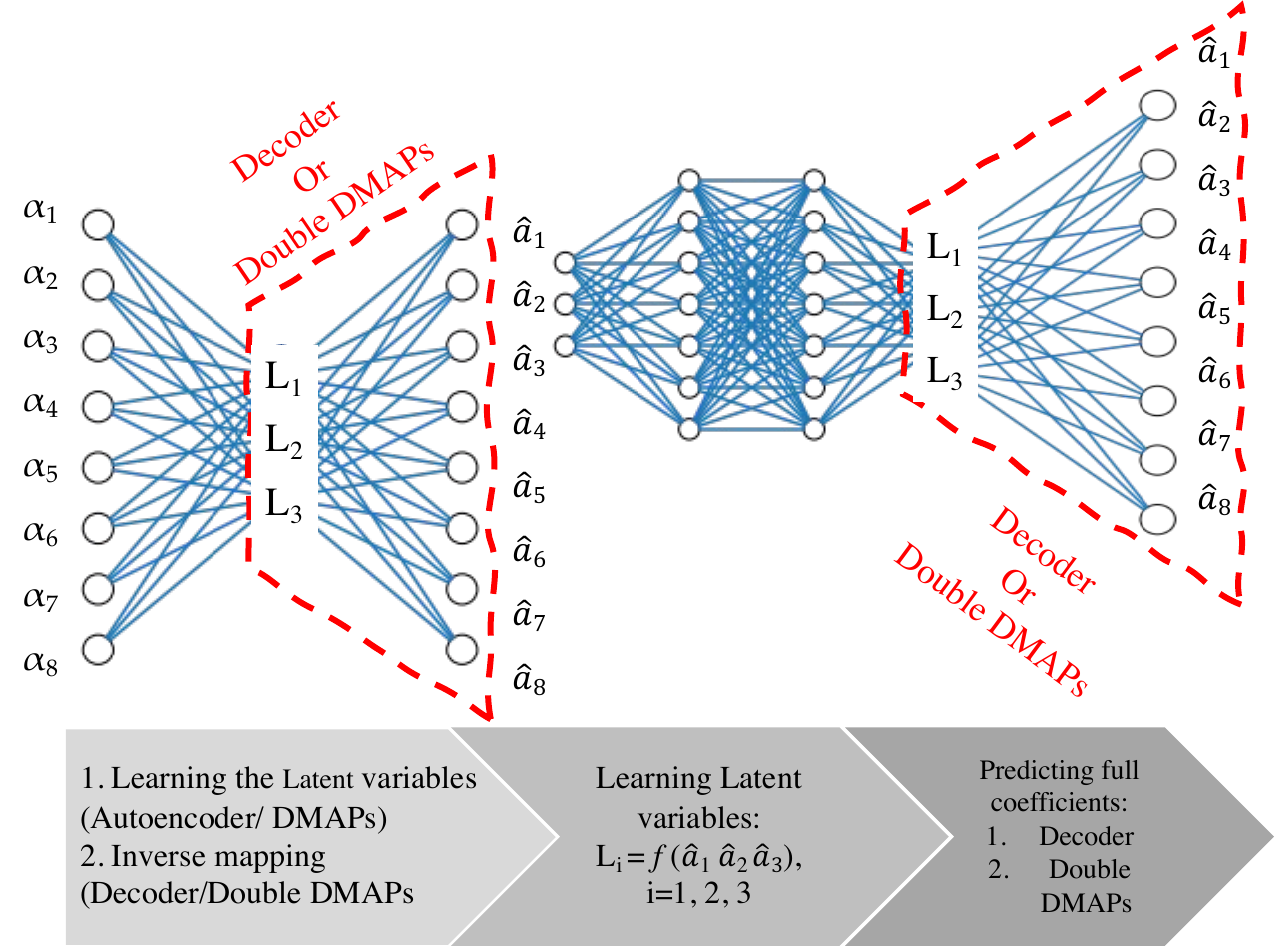}
    \caption{Schematic representation of computational workflow: First step includes learning the a minimum representation of the Approximate Inertial Manifold either with DMAPs or an autoencoder, as well as the inverse transformation, i.e. from the latent variables to the sine coefficients. Secondly, the latent variables are learned as a function of the leading three sines or POD coefficients, and finally, the full coefficients are predicted either with the decoder or Double DMAPS.}
    \label{fig:schema}
\end{figure}

The fact that the latent variables are one-to-one with the leading sine coefficients, allows us to recover the full sine coefficients in two distinct steps, schematically shown in Fig. \ref{fig:schema}. The first step is training the autoencoder. In the second step, we learn to infer the latent variables $\vect{L}$ from the leading sine coefficients, using either a feedforward neural network or Double DMAPs. 

Alternatively, the decoder part of the autoencoder can be used to compute an \textit{inverse-map}. This inverse map utilizes the leading Fourier modes, $\tilde{\vect{\alpha}}$, in which the dynamics have evolved, and the trained decoder, to find the latent autoencoder variables that minimize the algebraic optimization constraint
\begin{equation}
\label{eq:Optimization_Autoencoder}
   \operatorname*{argmin}_{\vect{L}}
|| \tilde{\vect{\alpha}} - \text{decoder}(\vect{L}) ||_F.
\end{equation}
In Equation \eqref{eq:Optimization_Autoencoder} the latent autoencoder variables are denoted as $\vect{L}$, the leading Fourier modes as $\tilde{\vect{\alpha}}$. After solving the optimization problem in Equation \eqref{eq:Optimization_Autoencoder} the decoder can be used to recover all the Fourier modes given $\vect{L}$.

This second use case of the autoencoder allows us to map from the leading Fourier modes to the latent space and back to the full Fourier models without the need of constructing an additional regression scheme. Once the latent variables are predicted, the decoder of the autoencoder or the inverse transformation from the DMAP to the ambient coordinates, is used to approximate the full set of reconstructed coefficients.
 
\section{Results}
\label{sec:Results}

\begin{figure}[ht!]
    \centering
        \includegraphics[width=0.8\textwidth]{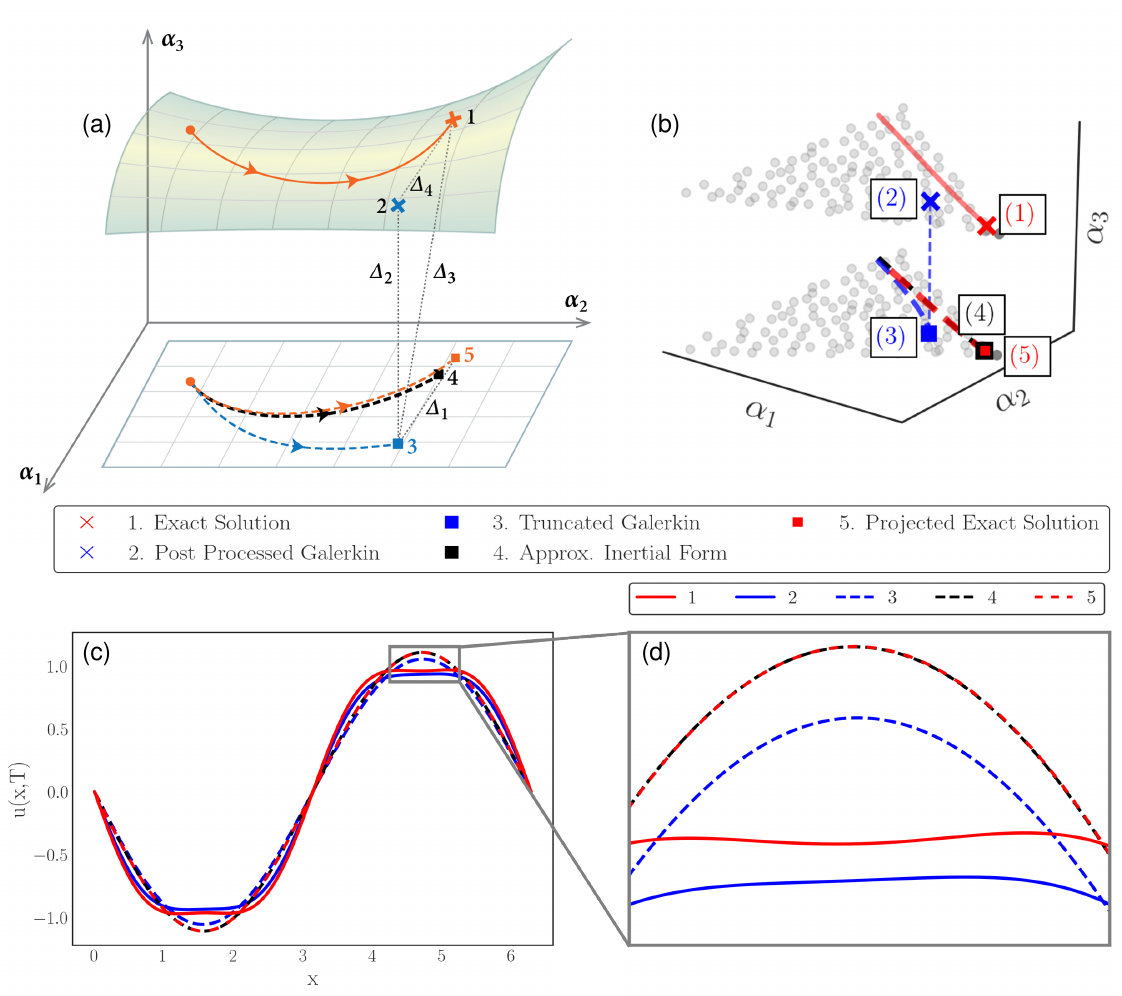}
	    \caption{(a) A schematic illustrating the benefits of the post-processing Galerkin methodology. A trajectory of the exact solution is shown \textit{on the manifold} in $a_1,a_2,a_3$ as a red solid trajectory, its final state denoted with an x marker and $(1)$. The projection of the exact solution in $a_1,a_2$ is shown with a red dashed line, its final state denoted as $(3)$, and a blue square. The trajectory integrated by using the approximate inertial form is shown with a black dashed line, and its final state is shown with a black square and denoted as $(4)$. The trajectory integrated by using the truncated Galerkin is shown with a blue dashed line, its final state denoted as $(3)$ and a blue square, and its post-processing (mapping) on the manifold, denoted as $(2)$, and indicated by a blue x marker. The dotted line $\Delta_1$ shows the distance between $(1)$ and $(5)$, the dotted line $\Delta_2$ shows the distance between $(2)$ and $(3)$, the dotted line $\Delta_3$ shows the distance between $(1)$ and $(4)$. The main \textit{premise} of post-processing Galerkin is that $\Delta_1$ and $\Delta_4$ are much smaller than $\Delta_2$. (b) The same components used in (a) are shown for the Chafee-Infante PDE. (c) The reconstructed solution in $u(x,T)$ for all possible options. (d) A blow-up of the reconstruction in $u(x,T)$.} 
	    \label{fig:TotalFigure}
\end{figure}

Before presenting our results we remind the reader, through the illustration in Fig. \ref{fig:TotalFigure} of the basic premise and the various errors associated with the post-processing Galerkin concept.  The main \textit{premise} is that the distance $\Delta_1$, between the projection of the true solution (point 5) and the truncated Galerkin solution (point 3) is much smaller than the distance $\Delta_3$ between point 3 and the true solution (point 1), the total error of truncated Galerkin approximation \citep{garcia1999approximate}.
This motivates the need for post-processing, which establishes that the distance $\Delta_4$ between point 1 and the post-processed Galerkin (point 2) is also much smaller than the total error (and comparable to $\Delta_1$) as shown in Fig. \ref{fig:TotalFigure}.



\subsection{{Euler-Galerkin vs. neural-network AIMs: Chafee-Infante}}
\label{sec:Chafee_Infante}
\begin{figure}[ht!]
    \centering
\begin{subfigure}[b]{0.4\textwidth}
         \centering
         \includegraphics[width=\textwidth]{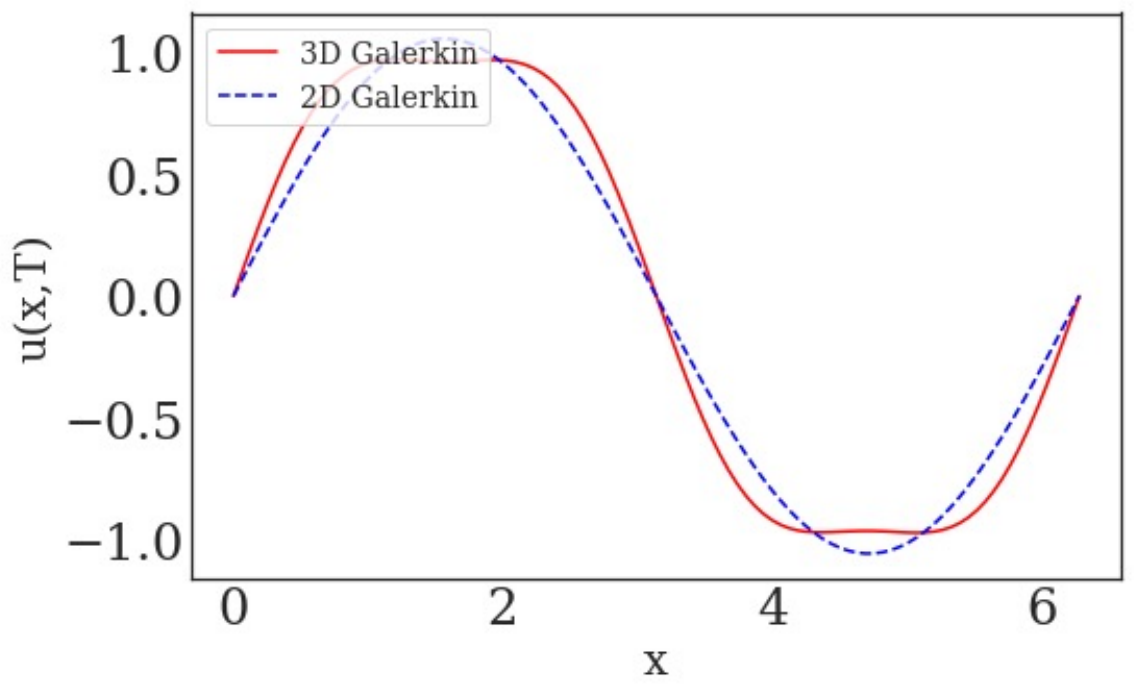}
         \caption{}
     \end{subfigure}
     \begin{subfigure}[b]{0.4\textwidth}
         \centering
         \includegraphics[width=\textwidth]{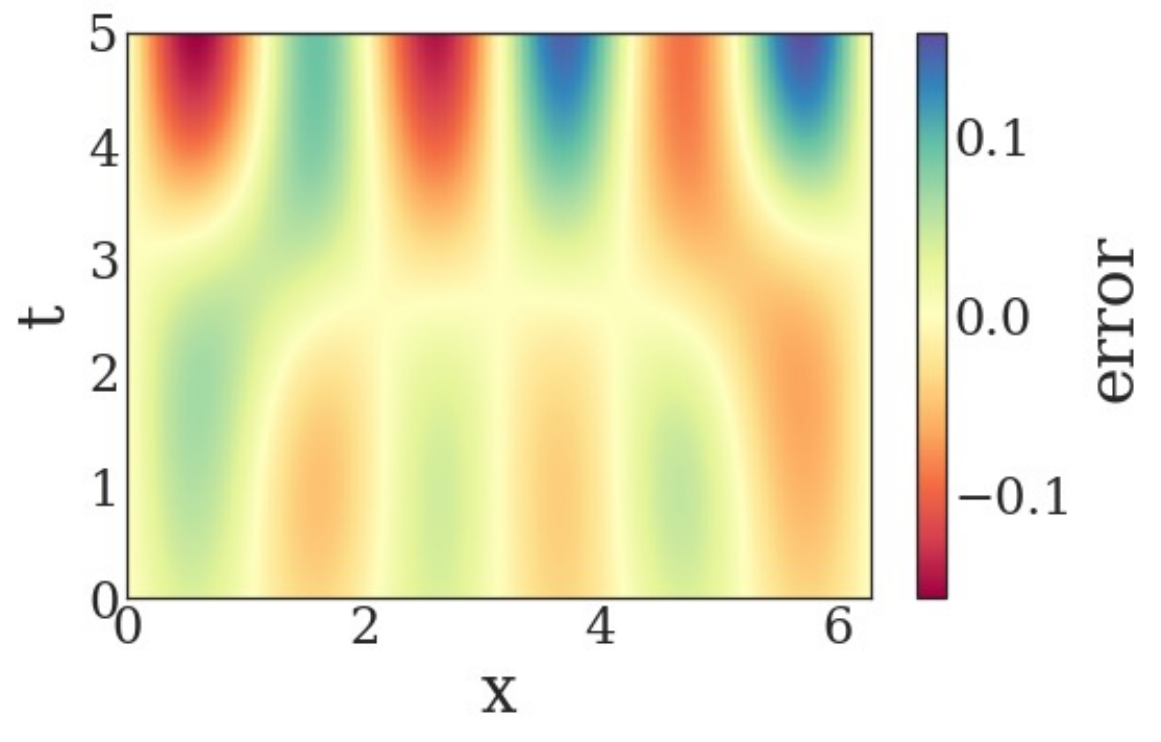}
         \caption{}
     \end{subfigure}
     \begin{subfigure}[b]{0.4\textwidth}
         \centering
         \includegraphics[width=\textwidth]{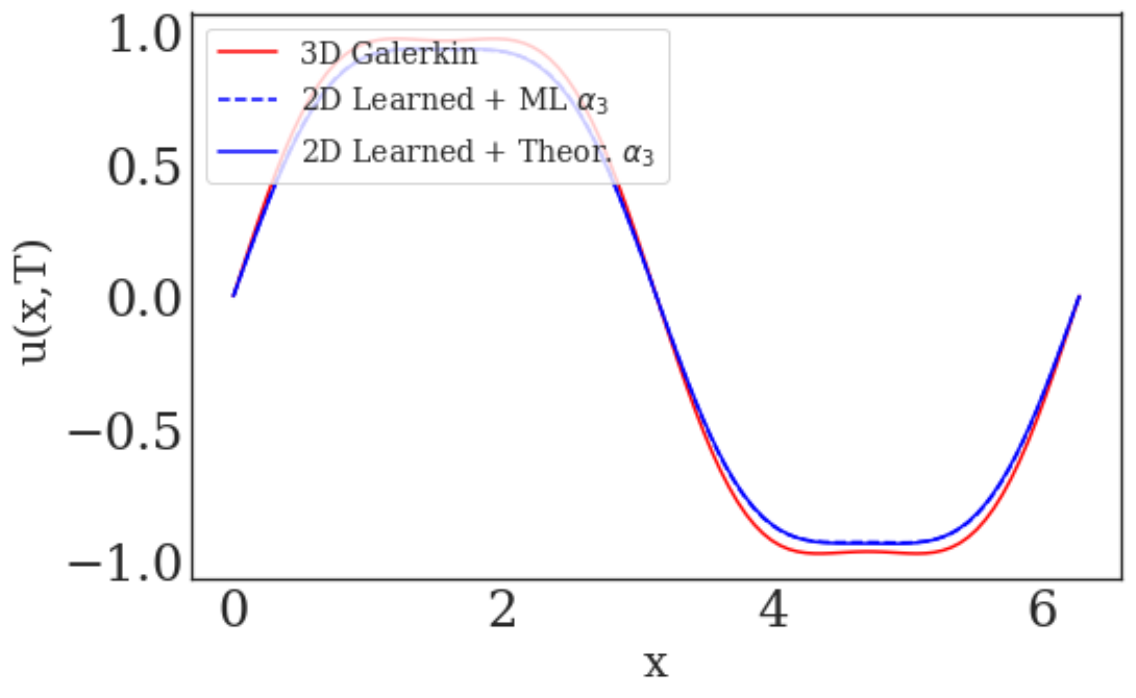}
         \caption{}
     \end{subfigure}
     \begin{subfigure}[b]{0.4\textwidth}
         \centering
         \includegraphics[width=\textwidth]{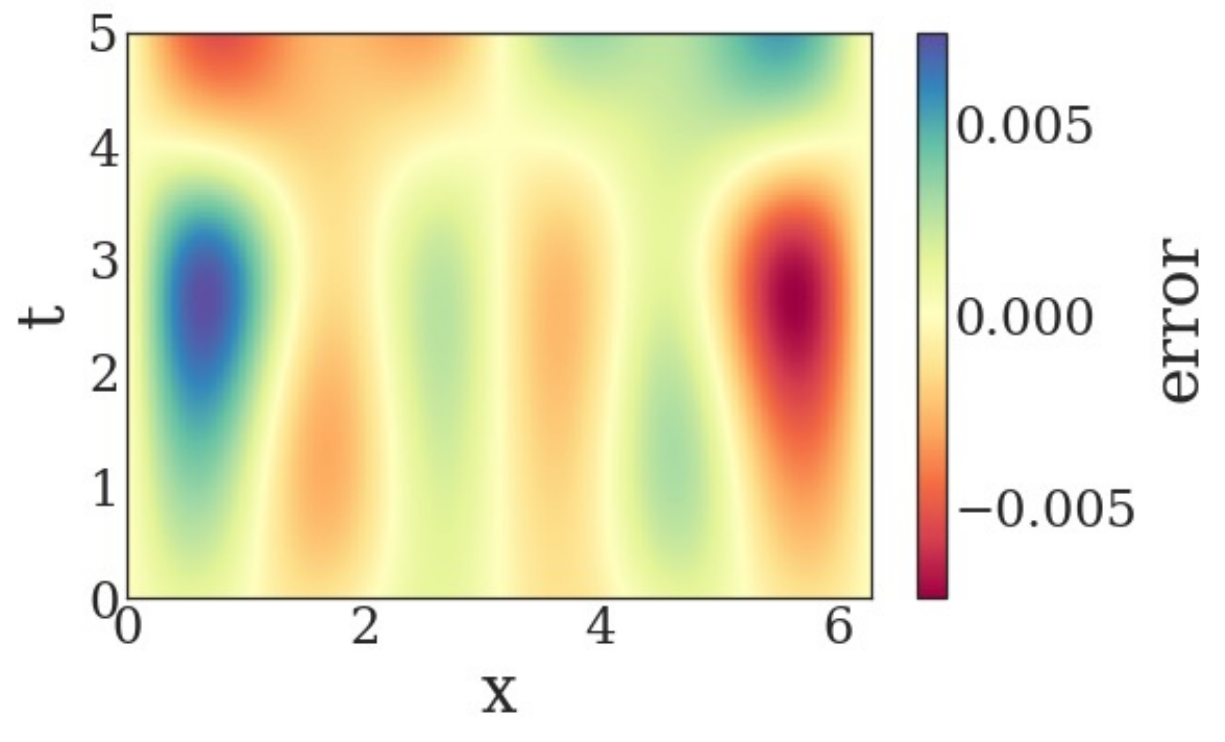}
         \caption{}
    \end{subfigure}
    \begin{subfigure}[b]{0.4\textwidth}
         \centering
         \includegraphics[width=\textwidth]{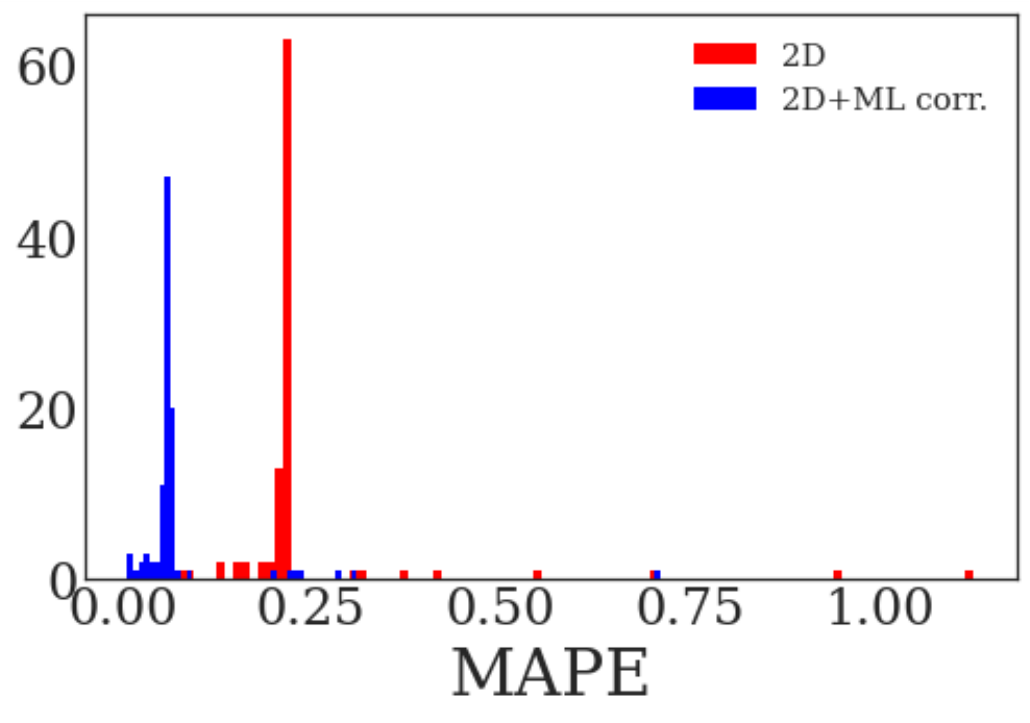}
         \caption{}
    \end{subfigure}
    \caption{Solution of Chafee Infante reconstructed in physical space; (a) The result of the 3D and the 2D Galerkin are shown in red and dashed black line respectively; (b) Percent error of reconstructed 2D Galerkin solution at each time-step; (c) Comparison of 3D Galerkin (red line) to the 2D Galerkin corrected with the neural network-derived term (dashed blue line); and the 2D Galerkin corrected with the theoretically derived $\alpha_3$; (d)  Percent error of reconstructed solution of the 2D ODE, corrected with the ML-derived $\alpha_3$, at each time-step; (e) Histogram of the mean absolute percent error of the 2D and ML-corrected 2D model at time $T=5$, for $100$ randomly selected initial conditions.}
    \label{fig:ChaInf2D3D}
\end{figure}

 For the Chafee-Infante we start by providing a comparison between the solution obtained with the three sine coefficients, here considered as the ground truth, (cf Fig.\ref{fig:ChaInf2D3D}a) and the truncated equations with the first two modes. The different post-processing schemes are applied to the solution of the truncated equations at the end of the desired integration. The comparison between the two is shown in Fig.\ref{fig:ChaInf2D3D} where  the reconstructed solution is shown with a blue dashed line and the ground truth simulation with a red line. The percent error along each step of the time integration until time $T=5$, is shown in Fig.\ref{fig:ChaInf2D3D}b.

 The solution of the 2D truncated dynamics is then corrected, using the value of $\alpha_3(T)$ as predicted by a neural network (described in Sec. \ref{subsec:LearnHarm}.1), using as inputs, the values of $\alpha_1(T)$ and $\alpha_2(T)$ at the final time-step, $t=T$. The results are shown in Fig.\ref{fig:ChaInf2D3D}c with a dashed blue line; included in the same figure, with a solid blue line, is the solution corrected with the theoretically (Euler-Galerkin AIM) derived value of $\alpha_3(T)$. The percent error along the integration time till time $T=5$ for the ML-derived $\alpha_3$ is shown in Fig.\ref{fig:ChaInf2D3D}d. Both, the ML-derived and the theoretical corrections help recover the accuracy and both lead to a mean absolute percent error (MAPE) of less than $1\%$. The mean absolute percent error (MAPE) is also computed at the same time instance ($T=5$) but for 100 randomly selected initial conditions, for the 2D and the ML-corrected 2D model. This is shown in Fig.\ref{fig:ChaInf2D3D}e, where the favorable effect of the correction on the mean absolute percentage error is clearly visualized. In these and in subsequent results, the MAPE refers to point-wise average of the absolute percentage error in each sample.
 
As an alternative, it is also possible to correct the \textit{learned} ODE in two dimensions, derived as described in \ref{sec:NNAIMandAIF}. The accuracy achieved is similar to the accuracy of the \textit{true} truncated 2D model.

\subsection{Kuramoto-Sivashinsky: Data-driven latent spaces and their AIMs}
The KS equation is selected in order to explore the application of the proposed methodology in cases where the minimum dimension of the Approximate Inertial Form (AIF) is not known a priori, although it has been argued to be three-dimensional \citep{jolly1990approximate}. Nevertheless, the truncated dynamics \textit{are not always quantitavely close to the actual behavior}. The latter will be addressed with ``Gray-Box" modeling, whereas the former is an important challenge in the implementation of post-processing Galerkin methods and will be addressed here with two different approaches: nonlinear manifold learning and in particular Diffusion Maps discussed in \ref{sec:dmaps}, and autoencoders discussed in \ref{sec:autoencoders}. 

\subsubsection{Learning the dimensionality of the latent space}

\noindent{\textbf{Autoencoders}}

A collection of data is sampled for the KS parameter value $\nu=33$, in various time instances of time-integration sufficiently close or on the global attractor. The data are used as inputs to an autoencoder and are reduced by the encoder into
 into a low dimensional bottleneck layer which parametrizes an approximation of the inertial manifold. It is then possible to map to the approximation of the high dimensional variables with the decoder. The encoder/decoder components of the network can be used independently as it will be demonstrated in a subsequent section to improve the accuracy of the reduced order model. 

\begin{figure}[ht!]
    \centering
    \includegraphics[width=\textwidth]{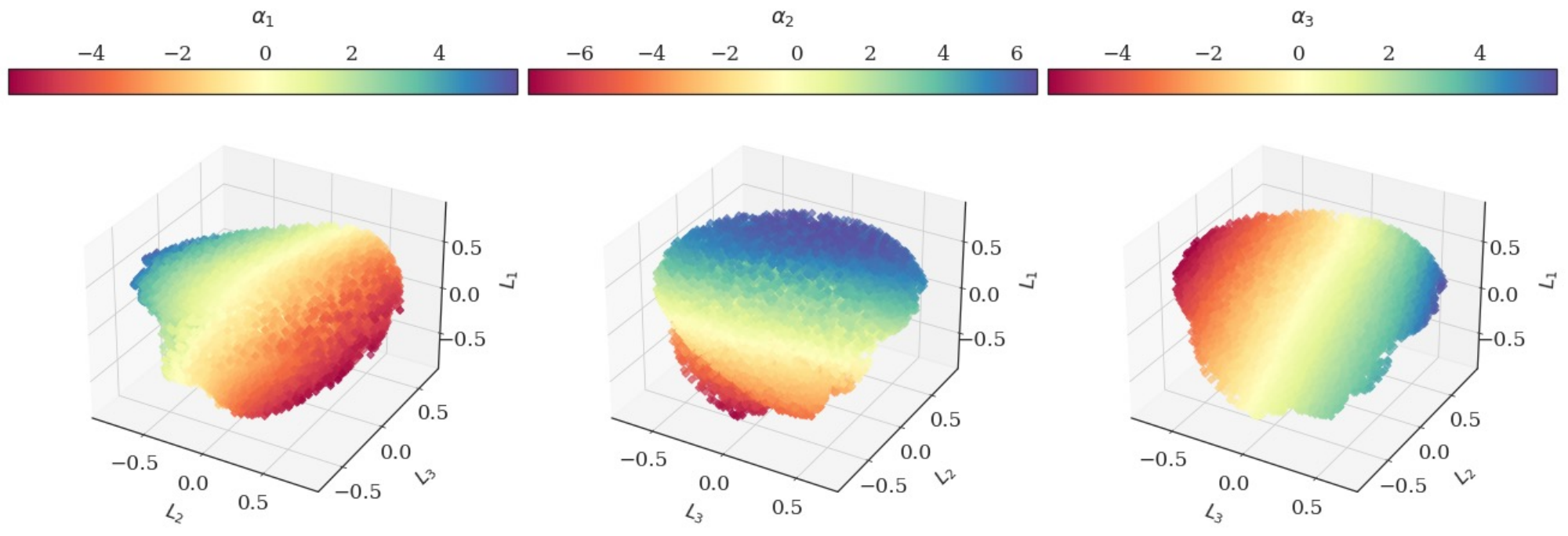}
    \caption{Latent variables of the autoencoder bottleneck layer; The three latent variables colored by the value of the first three sine coefficients, $\alpha_1$ (left), $\alpha_2$ (center) and $\alpha_3$ (right). Smoothness in color gradation suggests a one-to-one relation.} 
    \label{fig:autolatent}
\end{figure}

The three latent variables of the bottleneck layer are one-to-one functions of the first three sine coefficients, $\alpha_1, \alpha_2$ and $\alpha_3$. This is shown in Fig.\ref{fig:autolatent}, where the three bottleneck variables are plotted and colored according to the three sine coefficients. The smooth color variation suggests a one-to-one correlation between the latent and the ambient variables. It so happens that each one of the sine coefficients is one-to-one with each of the Diffusion Maps coordinates (the comparison is shown in the SI).

\begin{figure}[ht]
    \centering
    \includegraphics[width=7cm]{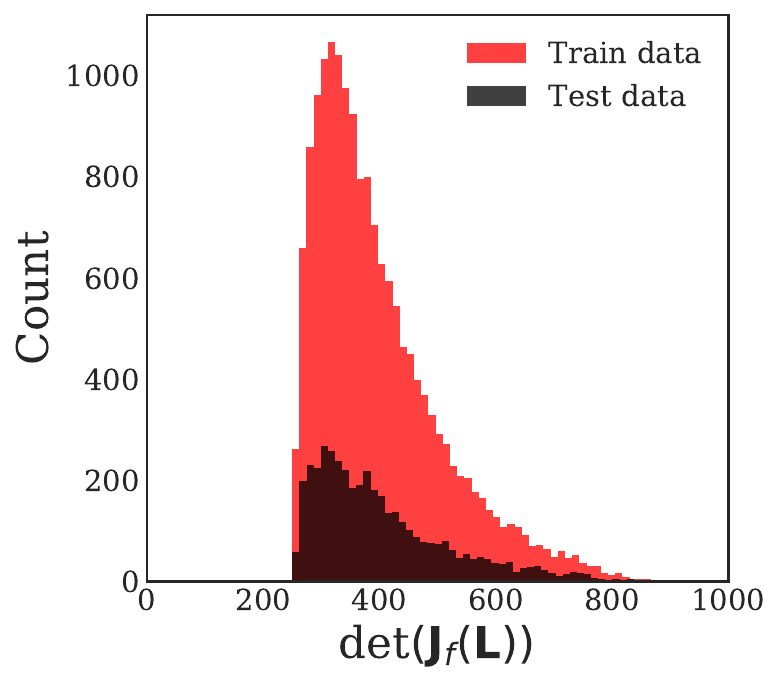}
    \caption{The histogram of the determinant of the Jacobian $\text{det}(\mathbf{J}_f(\mathbf{L}))$ computed along the training and test sets with automatic differentiation of the decoder.}
    \label{fig:Inverse_Function_Theorem_Histogram}
\end{figure}

 In Figure \ref{fig:Inverse_Function_Theorem_Histogram}, the histogram of the Jacobian determinant's values, $\text{det}(\mathbf{J}_f(\mathbf{L}))$, along the training and test data shows that this quantity is always positive and thus the mapping $f:\vect{L} \to \vect{ \tilde{\alpha}}$ is locally invertible. 

The one-to-one relationship between the leading sine coefficients and the autoencoder's latent variables $\vect{L}$ facilitates the second use case of the autoencoder we discussed earlier. 
This second use case utilizes the decoder to solve an \textit{inverse-problem} and map the leading sine coefficients $\tilde{\vect{\alpha}}$ to the autoencoder's latent space. Since, we showed that $f:\tilde{\vect{L}} \to \vect{\tilde{\alpha}}$ 
is a locally invertible map we can use the trained decoder and estimate $\vect{L}$ given $\tilde{\vect{\alpha}}$ by solving the optimization problem described in Equation \ref{eq:Optimization_Autoencoder}. As initial conditions to solve the optimization problem randomly sampled points from the training set were used. After optimization, the decoder can be used to reconstruct the remaining sine coefficients and from those the solution in $u(x,t)$ space, from the obtained values in autoencoder's latent space. 

\begin{figure}[ht!]
    \centering
\begin{subfigure}[b]{0.5\textwidth}
         \centering
         \includegraphics[width=\textwidth]{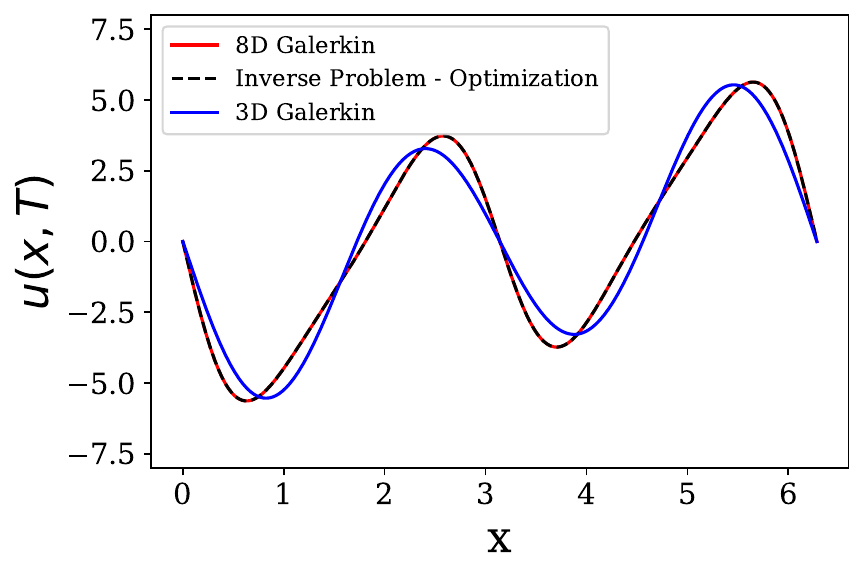}
         \caption{}
    \label{fig:optimization_trajectory}
     \end{subfigure}
     \begin{subfigure}[b]{0.32\textwidth}
         \centering
         \includegraphics[width=\textwidth]{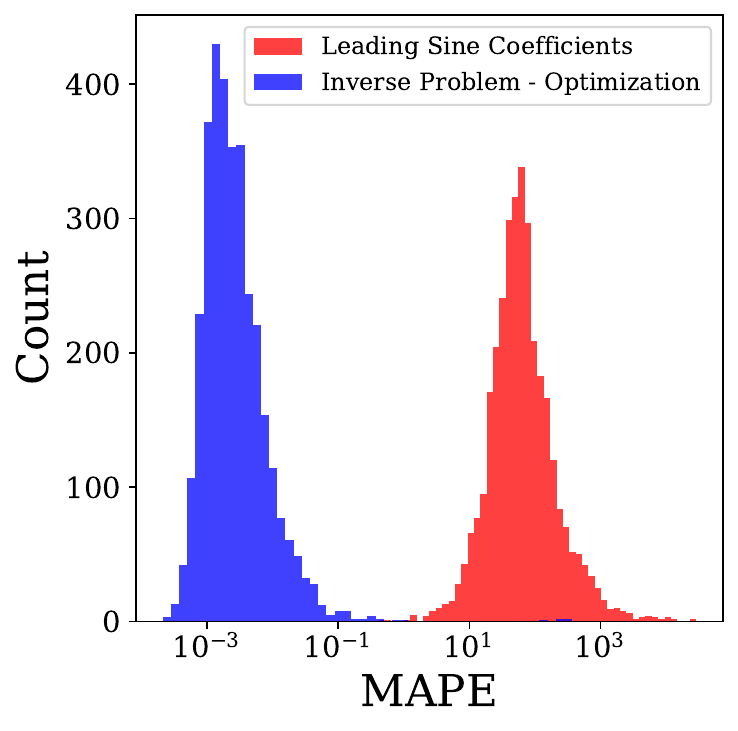}
         \caption{}
         \label{fig:y equals x}
     \end{subfigure}
    \caption{(a) A visual comparison for one data point between the true solution, the solution obtained after solving the inverse problem, and the truncated solution obtained by using the first three leading sine coefficients. (b) The mean absolute percent error (MAPE) across all the test points between the true solution and (i) the solution based on the leading sine coefficients (red histogram) (ii) the solution obtained after solving the inverse problem with optimization (blue histogram). }
    \label{fig:optimization_l2}
\end{figure}

In Figure \ref{fig:optimization_trajectory}, we contrast, for one reconstructed trajectory (i) the true solution $u(x,T)$ obtained from the full equations, (ii) the reconstructed solution based on the first three learned sine coefficients, and (iii) the reconstructed solution obtained by solving the inverse map and using the decoder to reconstruct the full solution. In Figure \ref{fig:optimization_l2}b the histograms show a comparison of the MAPE in $u(x,t)$ space between the true solution and the solutions based on (i) the three leading sine coefficients and (ii) the solution obtained after implementing the optimization step. \newline

\noindent{\textbf{{Diffusion Maps and their data-driven AIMs}}}

Diffusion Maps is implemented to encode the high dimensional data to a low dimensional manifold parametrized by three Diffusion Maps coordinates shown in Fig. \ref{fig:ks_dmaps} of the Appendix. The Diffusion Maps coordinates  $\phi_1$, $\phi_2$ and $\phi_3$ are one-to-one with the coefficients of the first three sine terms. This is shown in Fig. \ref{fig:ks_dmaps}, of the Appendix, by the smooth color transition in the diffusion maps plot when colored by $\alpha_1$, $\alpha_2$ and $\alpha_3$.


The sine coefficients, $\alpha_i$, are reconstructed with the help of Double DMAPS and by the decoder of the autoencoder as discussed in Secs. \ref{sec:dmaps} and \ref{sec:autoencoders} respectively. The MSE for the Double DMAPs approach is 0.00492, whereas for the autoencoder it is 0.0155. The precision of the autoencoder decreases for higher harmonics, which leads to the overall drop in accuracy of the reconstruction (this comparison is shown in the Sec. \ref{sec:SI}). Double DMAPs predicts accurately \textit{all} the coefficients (this comparison is shown in the Sec. \ref{sec:SI}).

\begin{figure}[ht!]
    \centering
    \begin{subfigure}[b]{0.45\textwidth}
         \centering
         \includegraphics[width=\textwidth]{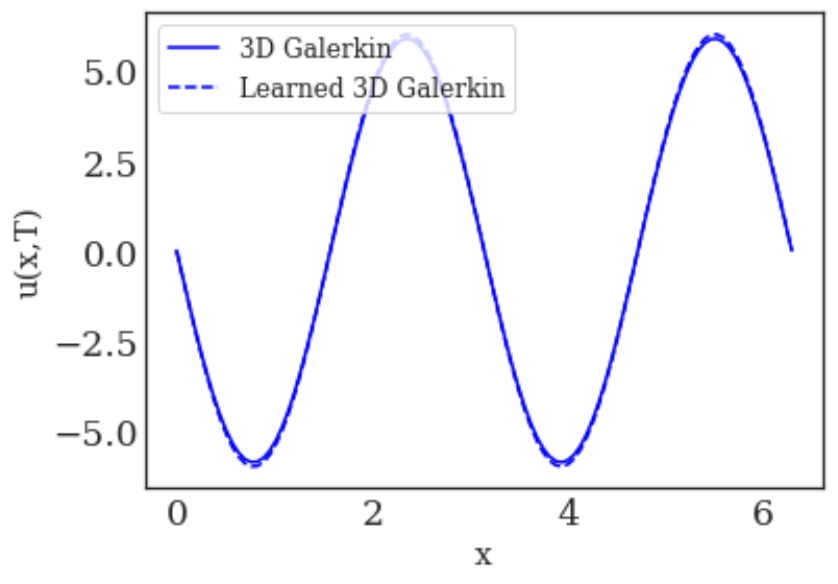}
         \caption{}
    \end{subfigure}
    \begin{subfigure}[b]{0.45\textwidth}
         \centering
         \includegraphics[width=\textwidth]{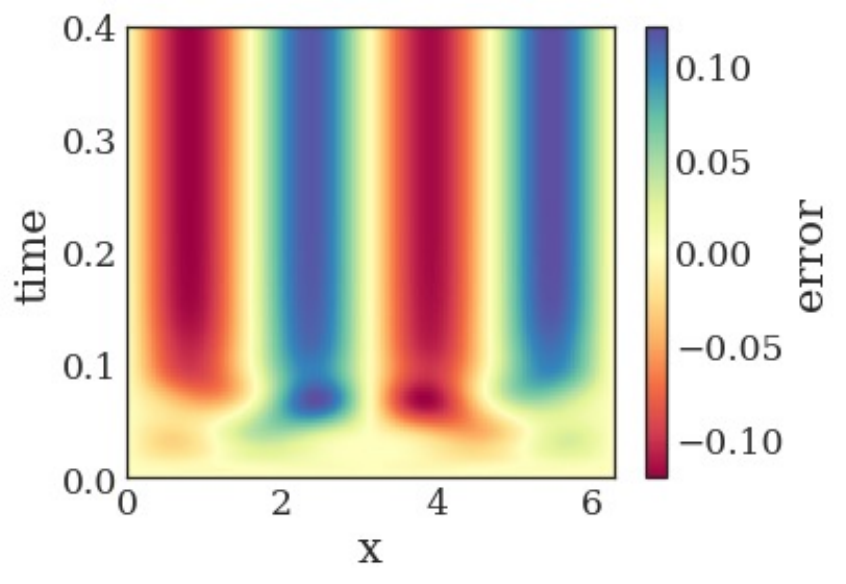}
         \caption{}
    \end{subfigure}
    \begin{subfigure}[b]{0.45\textwidth}
         \centering
         \includegraphics[width=\textwidth]{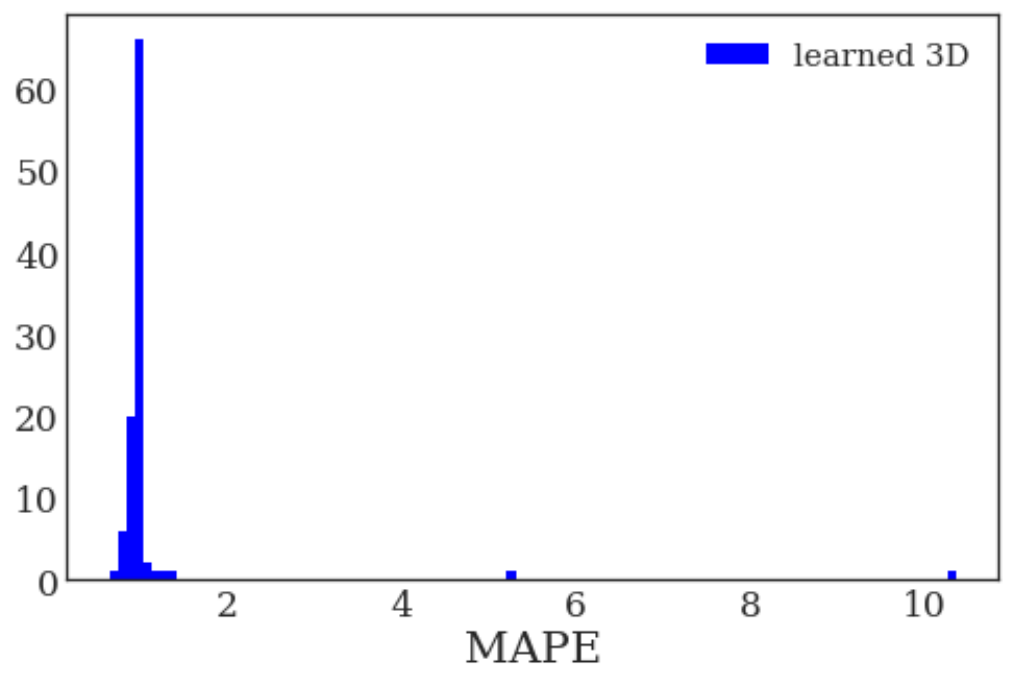}
         \caption{}
    \end{subfigure}
    \caption{(a) Reconstructed solution of the \textit{learned} 3D equation (broken line) and the actual 3D dynamics (solid line), showing almost perfect agreement.  (b) Percent error of the learned and true 3D solutions along the integration time. (c) Histogram of the mean absolute percent error of the learned 3D model at time $T$=0.06, for $100$ randomly selected initial conditions.}
\label{fig:auto3Dlearnedvs3D}
\end{figure}

\ek{\subsection{Data-driven post-processing Galerkin}}

\ek{Having established that the first three sine coefficients are one-to-one with the data-driven latent variables, the next step is to learn a data-driven ODE of the time evolution of the first three reconstructed sine coefficients as described in Sec. \ref{sec:NNAIMandAIF}.}

The feedforward neural network is trained using as input the values of $\hat{\alpha}_1$, $\hat{\alpha}_2$ and $\hat{\alpha}_3$, that are reconstructed by the latent space learning methods, i.e. autoencoders and DMAPs (the results of latent space identification and reconstruction of the hight-dimensional variables is presented in the SI). 
The predicted time-derivatives, the right-hand-side of the learned ODE, for each one of the sine coefficients are pictured in Fig. \ref{fig:reconstructed_rhs} (of the SI) versus the actual values of that components time derivative, $\dot{\alpha}$. The top row shows the predicted right-hand-side from the three sine coefficients resulting from the autoencoder, with $MSE = 9.5$. The respective predictions from the Double DMAPs reconstruction are shown on the bottom row, with $MSE=2.2$.

\ek{The neural network-derived approximation is then used in conjunction with an ODE solver, such as the Runge-Kutta, in order to integrate in time. The outcome of integration is reconstructed in physical space and compared to the outcome of the ground truth integration (in 8D) and also the reconstructed solution using only the first three modes of the ground truth. This is shown in Fig. \ref{fig:auto3Dlearnedvs3D}, alongside the error between the learned 3D ODE and the actual 3 modes of the Galerkin approximation, which demonstrates that the learned ODE predicts accurately the \textit{low dimensional time evolution of the first three modes}.}

\begin{figure}[ht!]
    \centering
    \includegraphics[width=\textwidth]{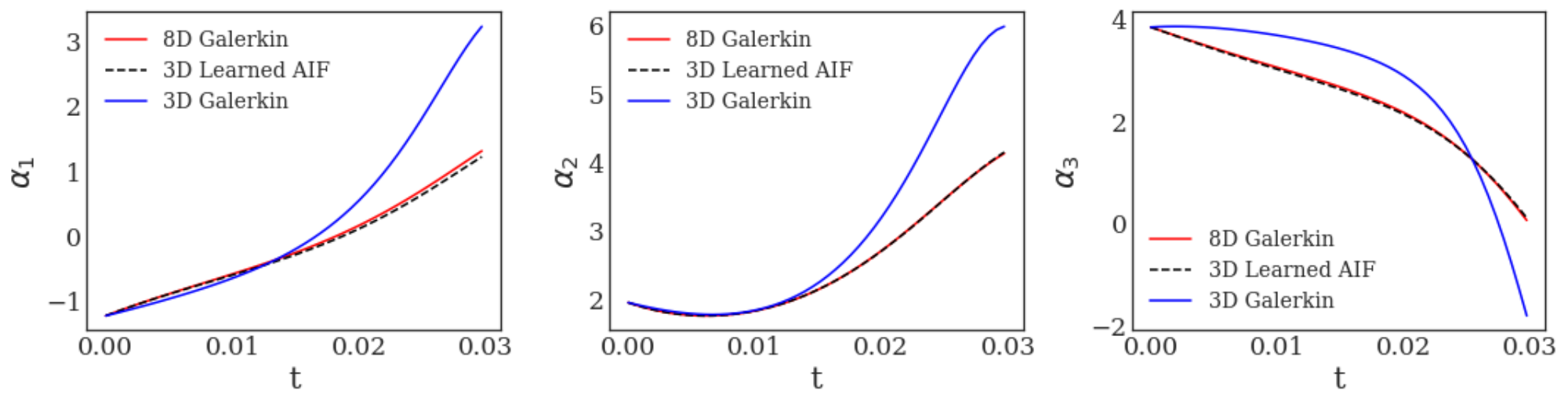}
    \caption{Left: Comparison of time evolution of the first three sine coefficients of the Galerkin discretization; 8-dimensional (red), learned 3-dimensional (black) the 3-dimensional (blue) Galerkin discretization}
    \label{fig:KSE8Dvs3Dcoeff}
\end{figure}

\begin{figure}
    \centering
    \begin{subfigure}[b]{0.45\textwidth}
         \centering
         \includegraphics[width=\textwidth]{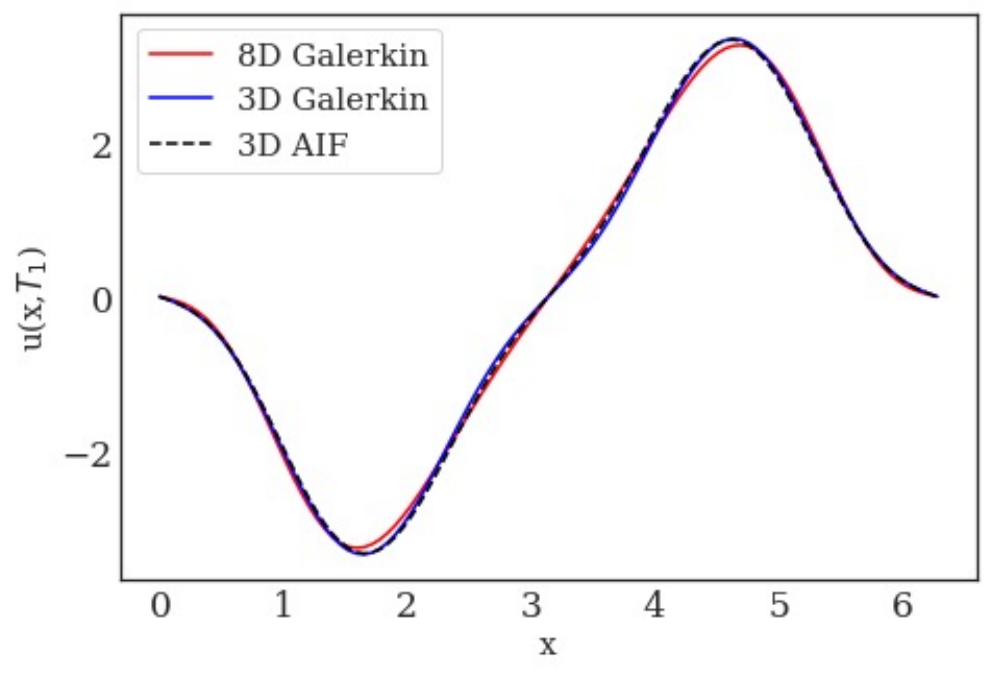}
         \caption{}
    \end{subfigure}
    \begin{subfigure}[b]{0.45\textwidth}
         \centering
         \includegraphics[width=\textwidth]{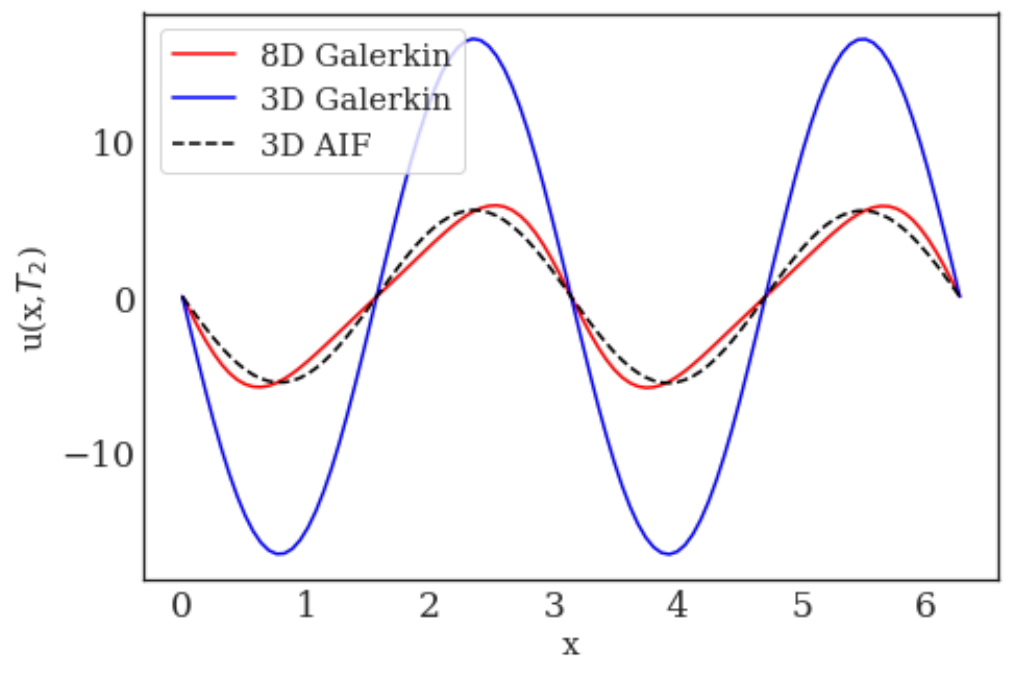}
         \caption{}
    \end{subfigure}
    \begin{subfigure}[b]{0.4\textwidth}
         \centering
         \includegraphics[width=\textwidth]{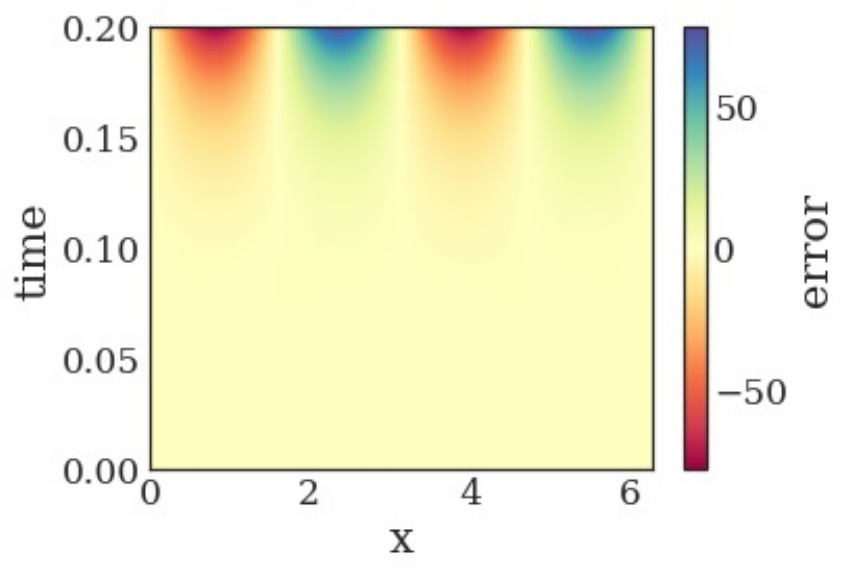}
         \caption{}
    \end{subfigure}
    \begin{subfigure}[b]{0.35\textwidth}
         \centering
         \includegraphics[width=\textwidth]{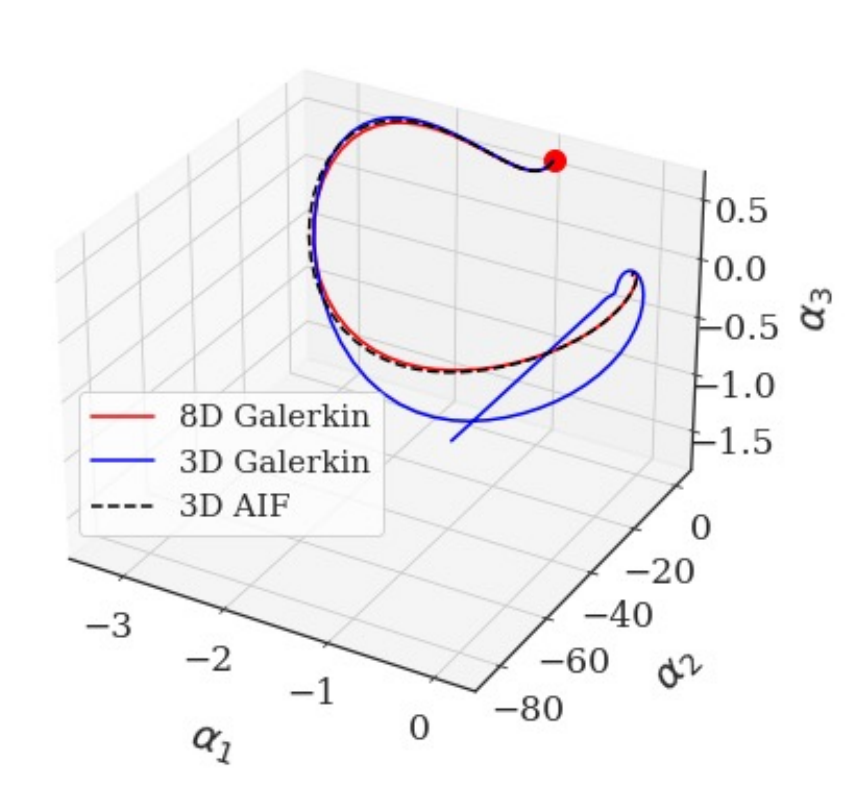}
         \caption{}
    \end{subfigure}

    \caption{ Comparison of the solution, in physical space, in two time instances (a) $T_1$, where the 3D and 8D dynamics are sufficiently close and (b) $T_2$, when they are quite far apart; (c) percent error between the truncated and the learned 3D AIF dynamics; (d) Comparison of time evolution of the first three sine coefficients of the Galerkin discretization : 8-dimensional (solid line), learned 3-dimensional (broken line) the 3-dimensional AIF (dotted line) Galerkin discretization}
\label{fig:KStruncomp}
\end{figure}

\subsubsection{When post-processing Galerkin works, when it does not work, and how to fix it}
\label{subsub:trunc}
It is worth looking into the time evolution of the first three modes, $\alpha_1$, $\alpha_2$ and $\alpha_3$ that result from the truncated 3D dynamics 
and compare it to the evolution of the first three terms of the full 8D dynamics and those of the learned AIF ODE. This is shown in Fig. \ref{fig:KSE8Dvs3Dcoeff}, where it becomes evident that the first three terms of the learned 3D ODE are close to the trajectory of the first three terms of the 8D Galerkin. In contrast, the truncated 3D dynamics deviate significantly,  past a certain point in time, from the ground truth dynamics. This observation suggests that using a post-processing scheme \textit{directly on the truncated equations} won't be able to correct the dynamics. This motivates us to use an ML correction so-called ``Gray-Box" model discussed further in the next section.

In essence, the post-processing Galerkin method relies on the premise that the solution of the truncated problem is reasonably close to the projection of the ground truth solution. Here it is demonstrated that even though the AIM is indeed three-dimensional, in the range of physical parameters examined, the truncated long term dynamics in three dimensional space are not accurate. This is demonstrated further in Fig. \ref{fig:KStruncomp}, where the solution in physical space is reconstructed from the 8D Galerkin, the 3D Galerkin and the learned 3D AIF ODEs, in different instances along the same trajectory. At initial stages of the trajectory, the solutions of the three methods are reasonably close, as shown in Fig. \ref{fig:KStruncomp}a. Later in time, the truncated 3D solution is growing quantitatively further apart from the ground truth, whereas the learned 3D AIF ODE follows closely the 8D dynamics (Fig.\ref{fig:KStruncomp}b). This is made clear in Fig. \ref{fig:KStruncomp}c, where the percent error between the learned and the truncated 3D solution is plotted along the trajectory. This can also be observed in phase space shown in Fig.\ref{fig:KStruncomp}d, on the right. The red point corresponds to the initial condition. The values of sine coefficients initially evolve in a similar manner but eventually, the truncated 3D dynamics deviate. \newline

\begin{figure}[ht!]
    \centering
    \begin{subfigure}[b]{0.45\textwidth}
         \centering
         \includegraphics[width=\textwidth]{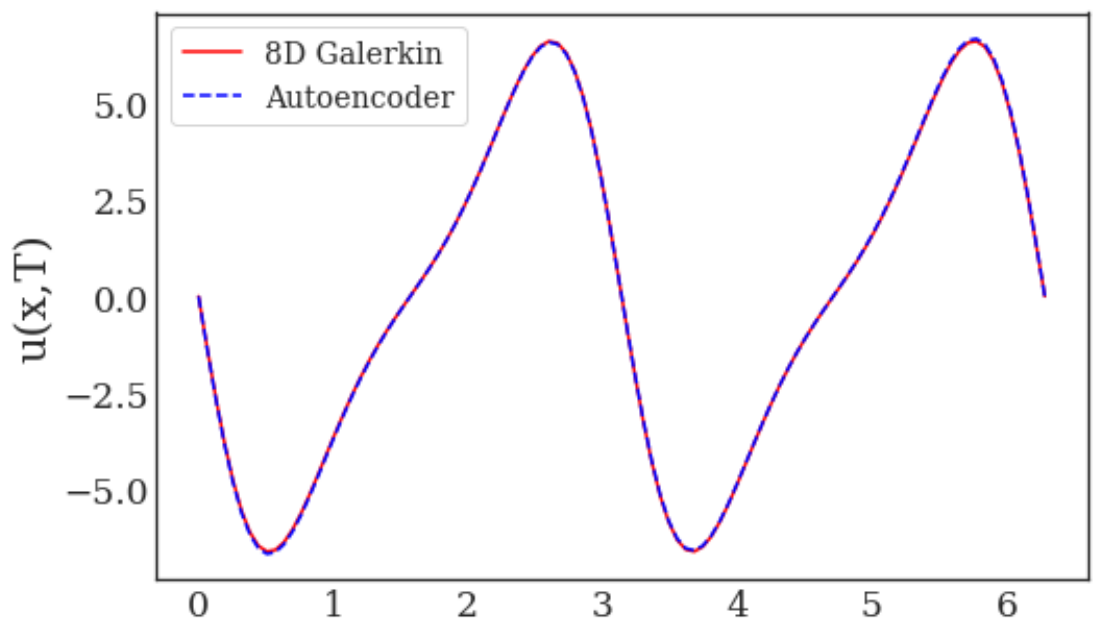}
         \caption{}
    \end{subfigure}
    \begin{subfigure}[b]{0.45\textwidth}
         \centering
         \includegraphics[width=\textwidth]{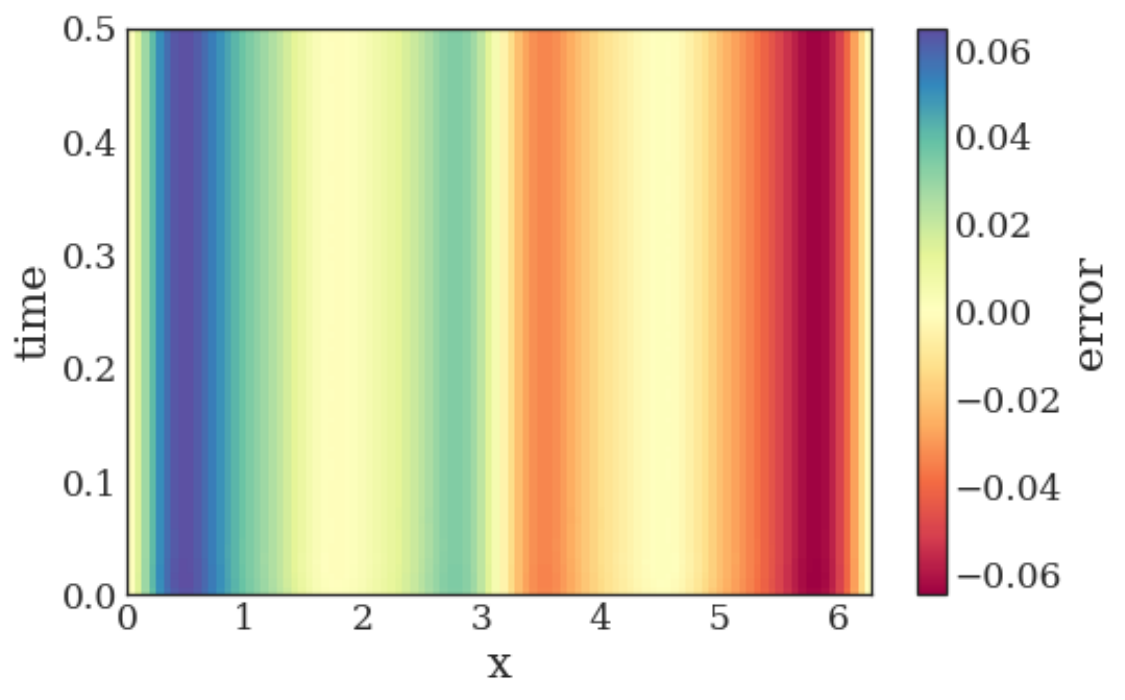}
         \caption{}
    \end{subfigure}
    \begin{subfigure}[b]{0.45\textwidth}
         \centering
         \includegraphics[width=\textwidth]{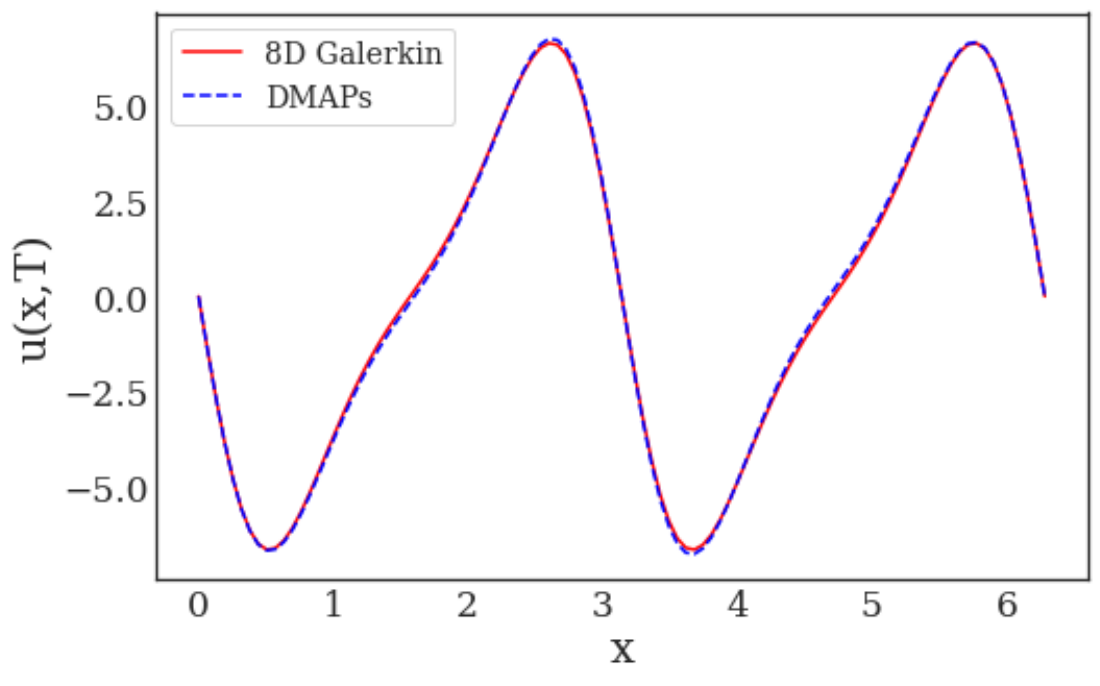}
         \caption{}
    \end{subfigure}
    \begin{subfigure}[b]{0.45\textwidth}
         \centering
         \includegraphics[width=\textwidth]{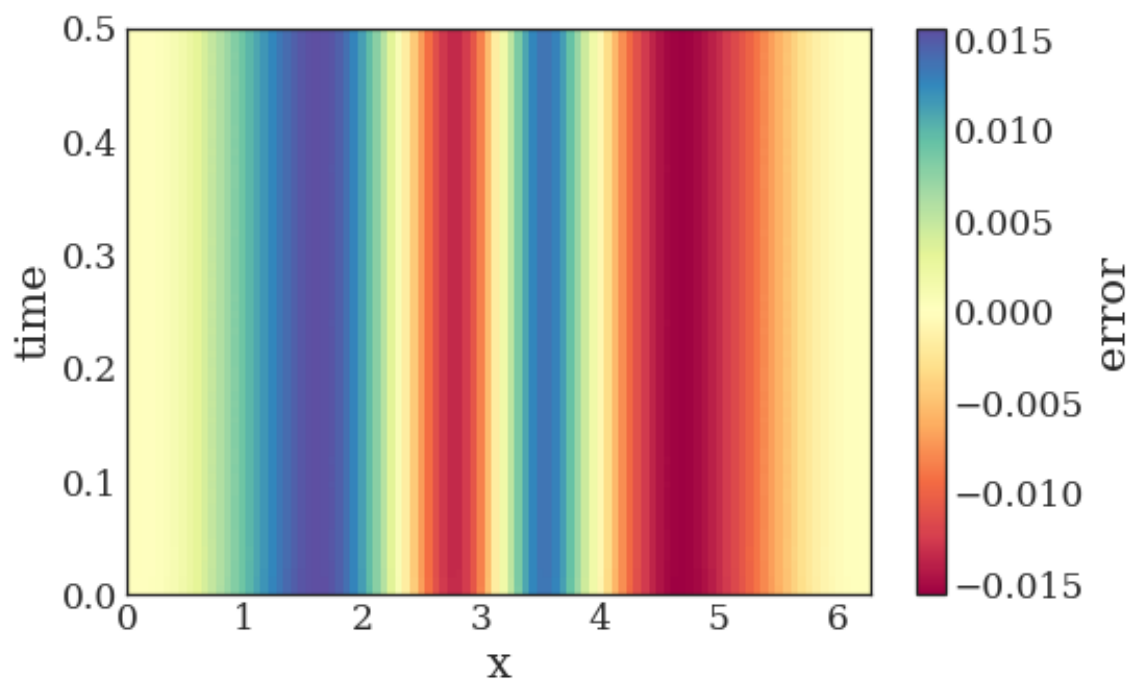}
         \caption{}
    \end{subfigure}
    \begin{subfigure}[b]{0.45\textwidth}
         \centering
         \includegraphics[width=\textwidth]{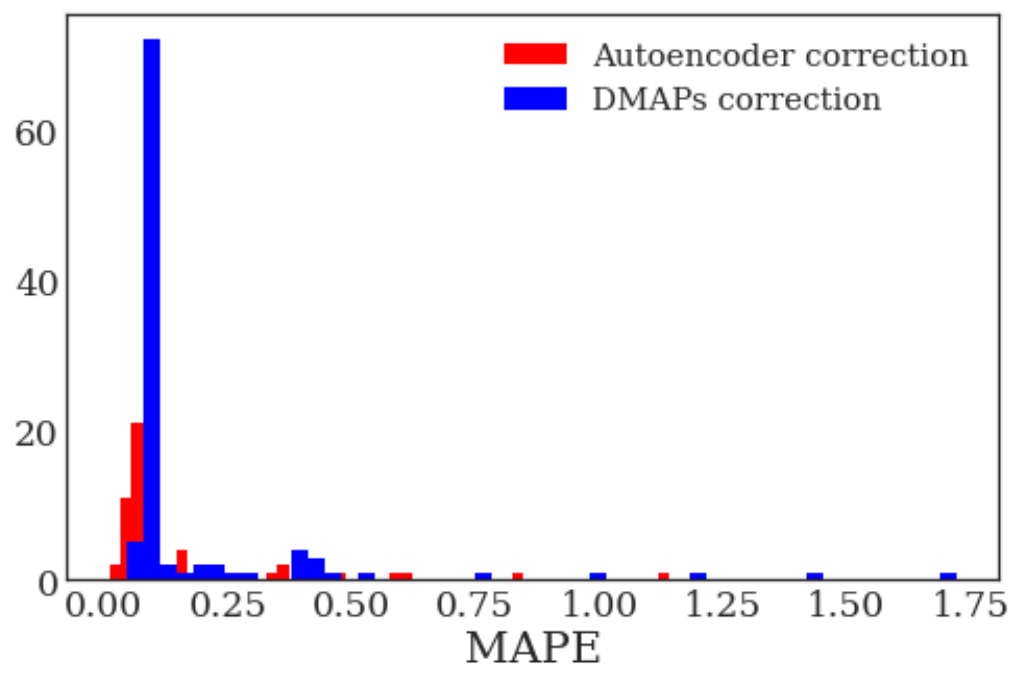}
         \caption{}
    \end{subfigure}
    \caption{(a) Comparison, at $\alpha$=33 and $T=0.5$ of the reconstructed Kuramoto-Sivashinsky solution, between the 8-dimensional Galerkin (solid red line) and the 3-dimensional learned AIF ODE corrected with the decoder-derived higher harmonics terms (broken blue line); (b) percent error over the time integration interval. (c) Comparison, at $\alpha$=33 and $T=0.5$ of the reconstructed Kuramoto-Sivashinsky solution, between the 8-dimensional Galerkin (solid red line) and the 3-dimensional learned AIF ODE corrected with the DMAPs-derived higher harmonics terms (broken blue line); (d) Percent error along the time-span of integration. (e) Histograms of MAPE of the autoencoder-corrected and DMAPs-corrected solution at time $T$=0.5, for $100$ randomly selected initial conditions.}
\label{fig:corrected}
\end{figure}

\begin{figure}[ht!]
    \centering
    \begin{subfigure}[b]{0.45\textwidth}
         \centering
         \includegraphics[width=\textwidth]{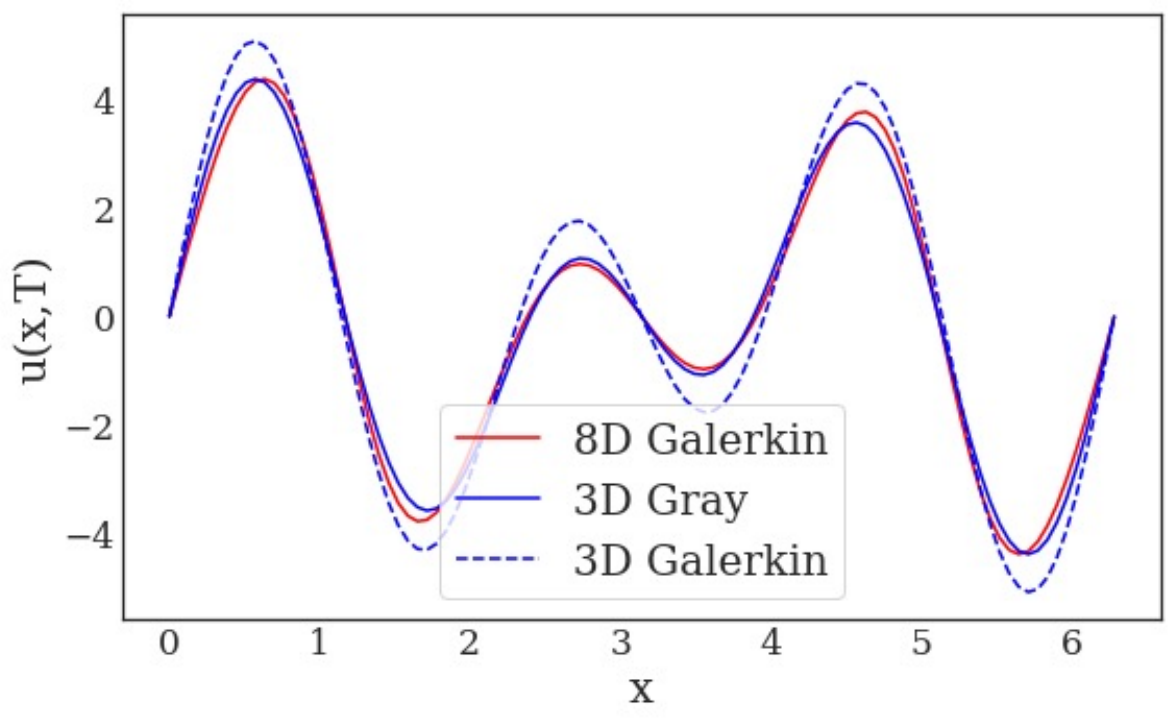}
         \caption{}
    \end{subfigure}
    \begin{subfigure}[b]{0.45\textwidth}
         \centering
         \includegraphics[width=\textwidth]{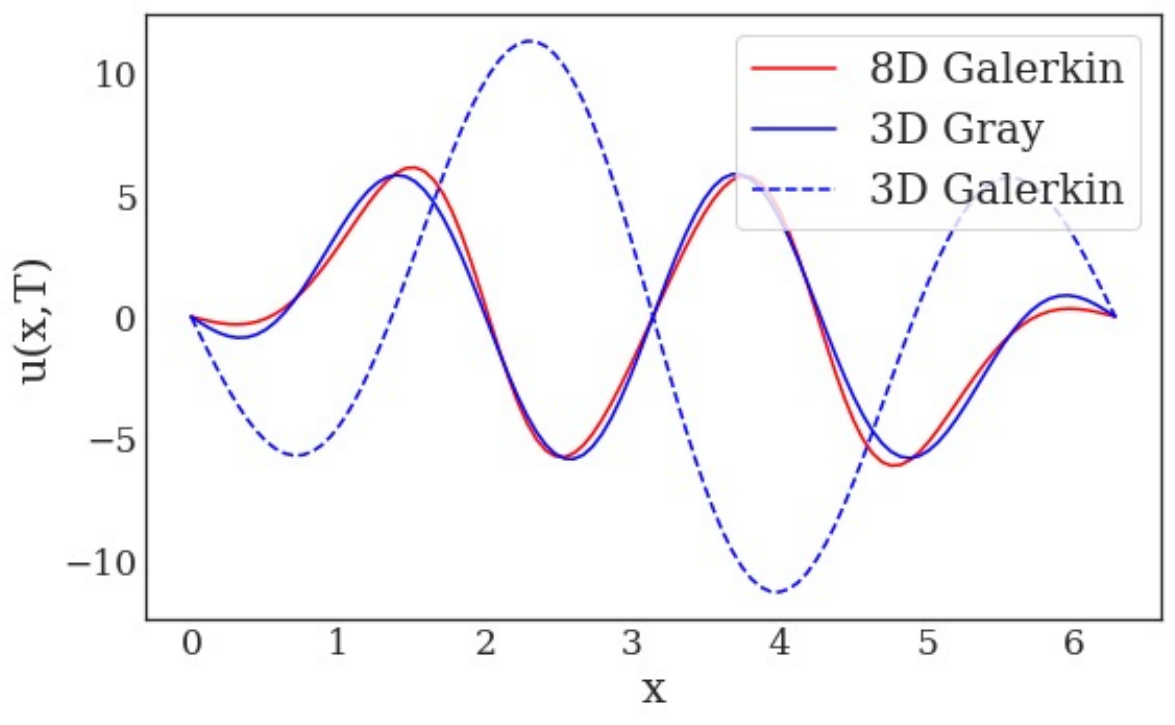}
         \caption{}
    \end{subfigure}
    \begin{subfigure}[b]{0.5\textwidth}
         \centering
         \includegraphics[width=\textwidth]{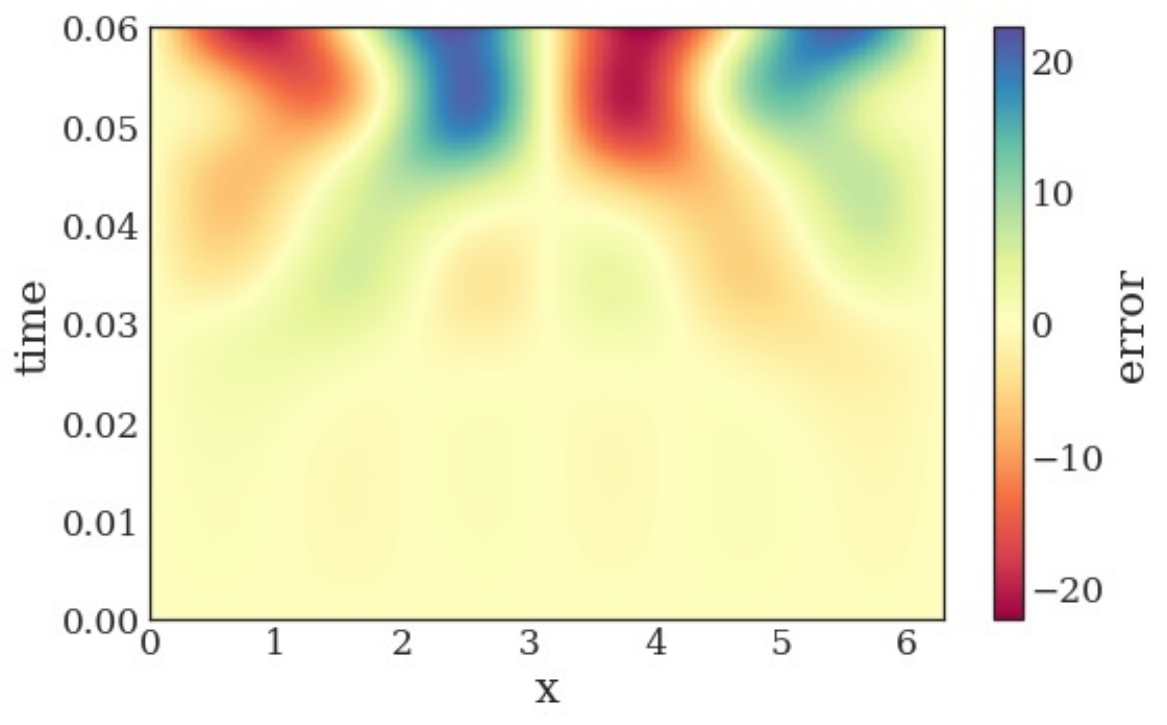}
         \caption{}
    \end{subfigure}
    \begin{subfigure}[b]{0.45\textwidth}
         \centering
         \includegraphics[width=\textwidth]{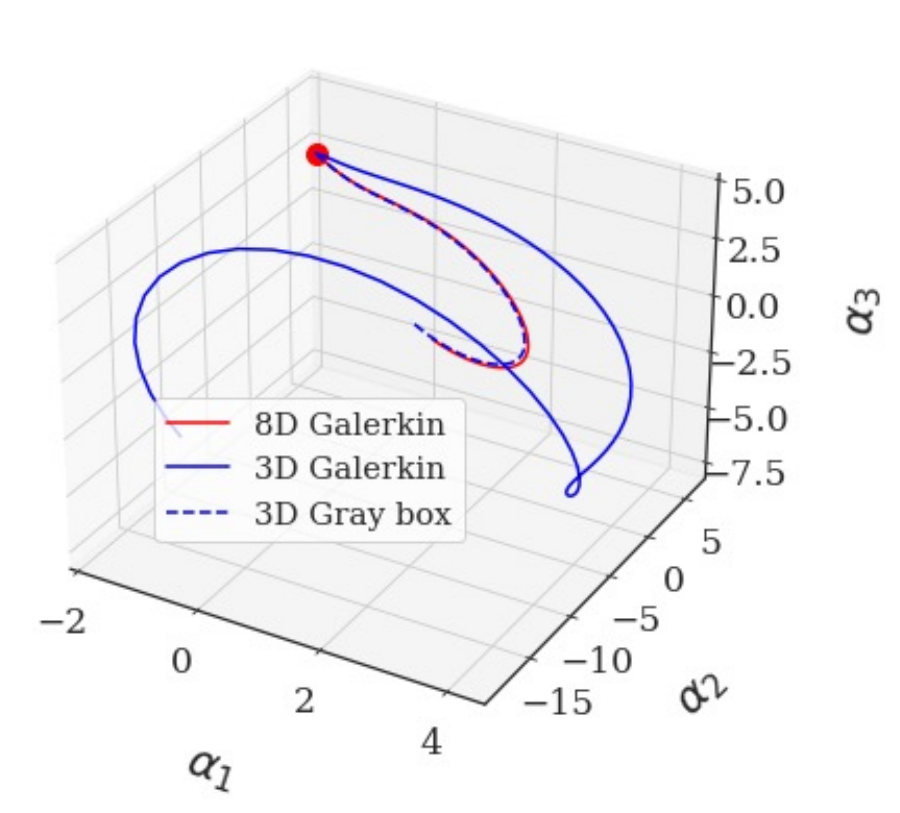}
         \caption{}
    \end{subfigure}
    \begin{subfigure}[b]{0.45\textwidth}
         \centering
         \includegraphics[width=\textwidth]{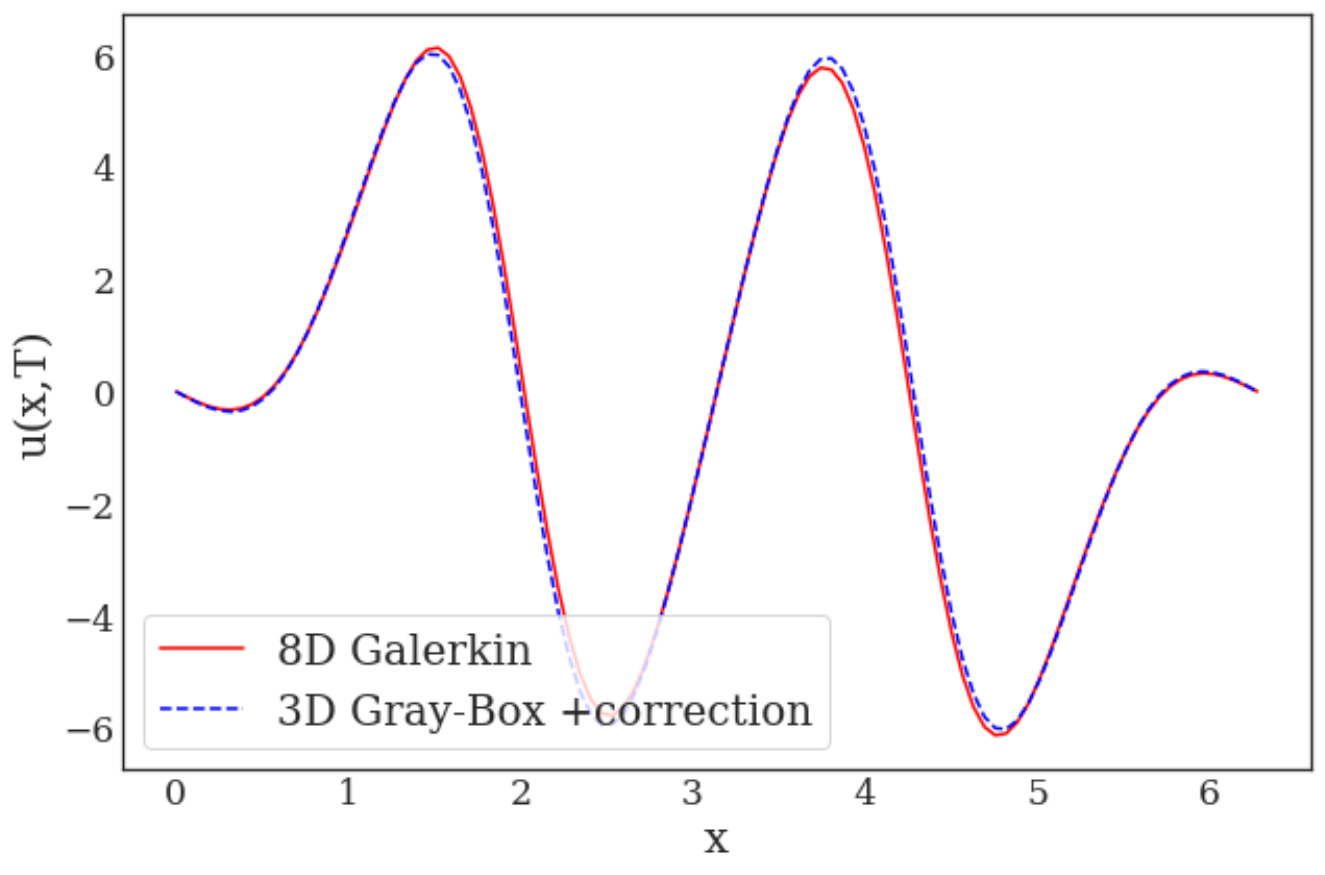}
         \caption{}
    \end{subfigure}

    \caption{Gray-Box correction of 3D Galerkin dynamics: Comparison of the solution, in physical space, in two time instances (a) $T_1=0.02$, where the 3D and the 8D Galerkin dynamics are sufficiently close and (b) $T_2=0.05$, when they are quite far apart; (c) Percent error between the truncated and the Gray-Box 3D dynamics; (d) Comparison of time evolution of the first three sine coefficients of the Galerkin discretization : 8-dimensional (solid red line), Gray-Box 3-dimensional (broken blue line) truncated 3-dimensional (blue solid line) Galerkin discretization (e) Comparison of the solution, in physical space, of the corrected Gray-Box with the 8D Galerkin at $T_2=0.05$}
\label{fig:gray}
\end{figure}

\noindent{\textbf{Data-driven post-processing Galerkin}}

To recover the values of all the $\alpha_i$s, necessary for accurate reconstruction, the first step involves predicting the latent variables, either the bottleneck variables from the autoencoder or alternatively the DMAPs latent coordinates. One way to achieve this, is with a feedforward neural network with three inputs (the the first three sine coefficients) and three outputs (the latent variables). In this implementation, the neural network consists of 5 hidden layers with 80 neurons each and a $tanh$ activation function, implemented in tensorflow \citep{tensorflow2015-whitepaper}. The mean squared error is used as the loss function along with the Adam optimizer.

It is then possible to employ either Double DMAPs, in the case of DMAPs, or the decoder of the autoencoder, and predict the corresponding $\alpha_i$s with $MSE=0.09$ for both cases. The reconstructed solution in physical space is compared to the ground truth in Fig. \ref{fig:corrected}.\hfill 

\newpage

\noindent{\textbf{When the PPG assumptions fail}}

If we have a accurate low-dimensional observation we can correct, in principle theoretically with Euler-Galerkin, or in practice with machine learning approaches as described above. If this is not available, then we proceed to improve the AIF itself though the Gray-Box approach.

The method's performance is demonstrated in Fig. \ref{fig:gray} for two cases. In the first case (cf Fig. \ref{fig:gray}a), the 3D, 8D, and corrected Gray-Box dynamics are shown for the reconstructed physical space solution, at $T_1=0.02$, when the truncated 3D dynamics are close to the ground truth. At a later time-step, at $T_2=0.05$ (cf Fig. \ref{fig:gray}b), the truncated dynamics have deviated far from the truth. The Gray-Box model corrects the deviation in both cases and accurately captures the ground truth with the addition of post-processing terms, as seen in  Fig. \ref{fig:gray}e. 

\subsection{Using POD coefficients to parametrize the IM/AIM}
\label{sec:PODcoef}
\begin{figure}[ht!]
    \centering
    \begin{subfigure}[b]{0.45\textwidth}
         \centering
         \includegraphics[width=\textwidth]{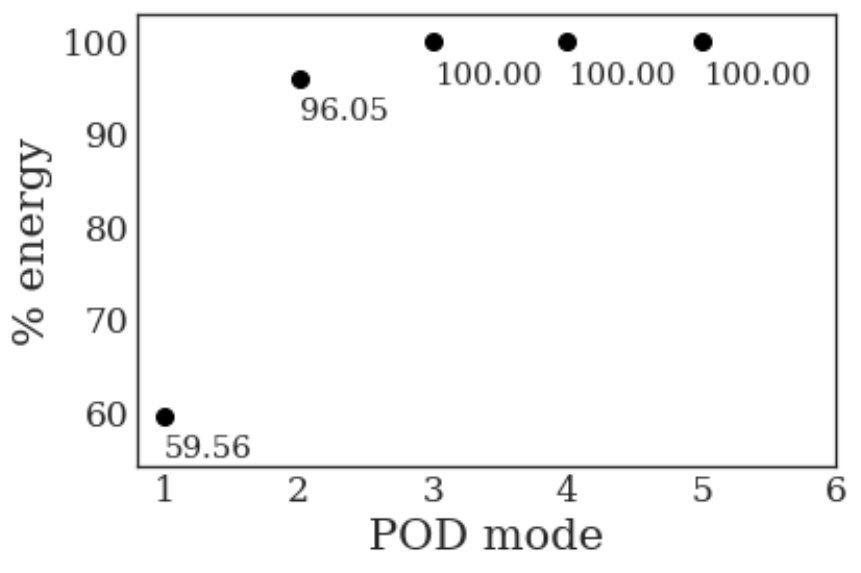}
         \caption{}
         \label{fig:energy}        
    \end{subfigure}
    \begin{subfigure}[b]{0.45\textwidth}
         \centering
         \includegraphics[width=\textwidth]{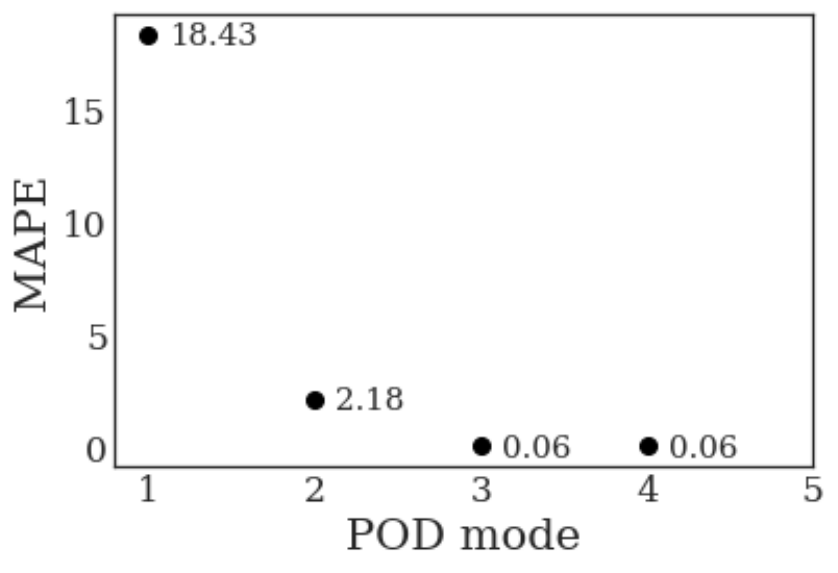}
         \caption{}
          \label{fig:mape}        
    \end{subfigure}
    \caption{(a) $\%$ energy of the data contained by progressively increasing POD basis size. 3 POD modes represent 99.99 $\%$ of the variance (b) MAPE of the dataset projected on progressively increasing POD basis, with respect to the original. When condidering 3 POD modes, the MAPE drops to $0.06 \%$.}
\end{figure}

Here, the implementation of the proposed workflow is presented in the case where the manifold is parametrized by data-driven POD coefficients, rather than sine coefficients, for the Chafee-Infante equation.
To start with, the POD modes that contain the greatest percentage of variance of an ensemble of solutions in physical space, are identified. Three POD modes represent $99.99 \%$ of the energy of the dataset (cf. Fig. \ref{fig:energy}), defined as the percentage of the cumulative sum of the leading three eigenvalues over the sum of all the eigenvalues. 

\begin{figure}[ht!]
    \centering
    \includegraphics[width=\textwidth]{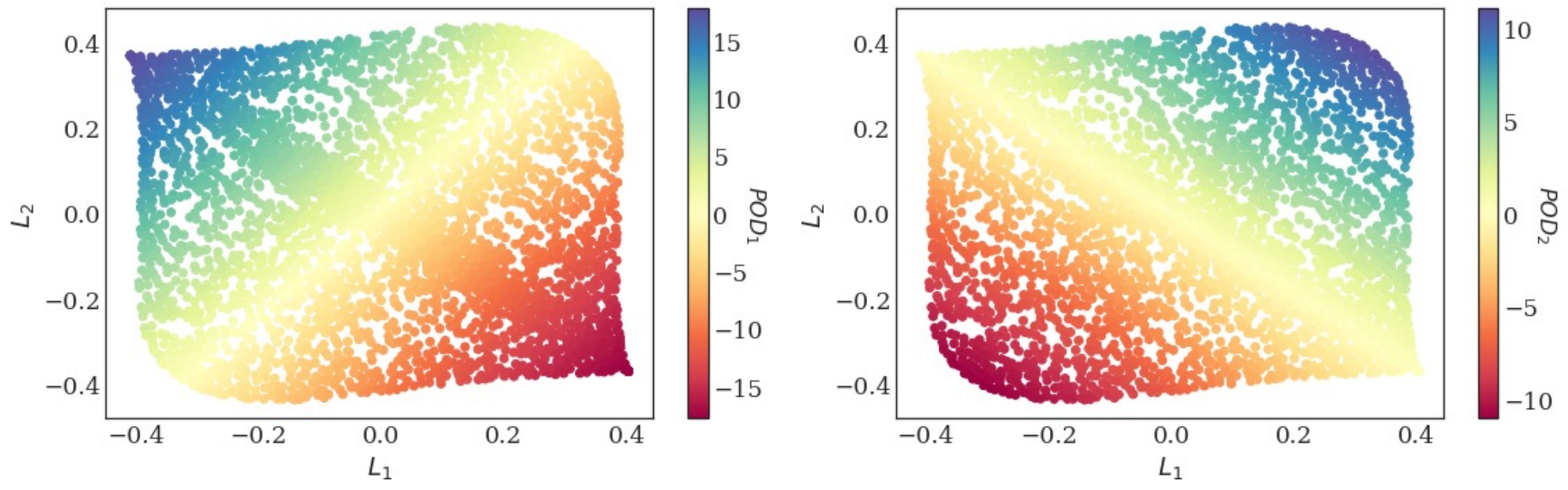}
    \caption{Latent variables discovered by the autoencoder, whose inputs are POD coefficients values. The smooth color transition implies that the latent variables are 1-to-1 with the leading POD coefficients.}
    \label{fig:POD1}
\end{figure}

The original dataset is then projected on the first three modes, leading to each solution vector, being represented by three coefficients. The mean absolute percent error of the dataset, projected on a basis consisting of 3 POD vectors, is $0.06 \%$ (cf. Fig.\ref{fig:mape}). We use this collection of POD coefficients, to discover the latent variables, here with an autoencoder with a 2-neuron bottleneck layer. The mean absolute percentage error achieved for the autoencoder-reconstructed POD coefficients is $1.2\%$. The latent variables are one-to-one with the two leading POD coefficients, as is evident in Fig. \ref{fig:POD1}, where they are plotted and coloured according to the values of the coefficients. The smooth colour transition is indicative of the one-to-one relationship.

\begin{figure}[ht!]
    \centering
    \includegraphics[width=0.5\textwidth]{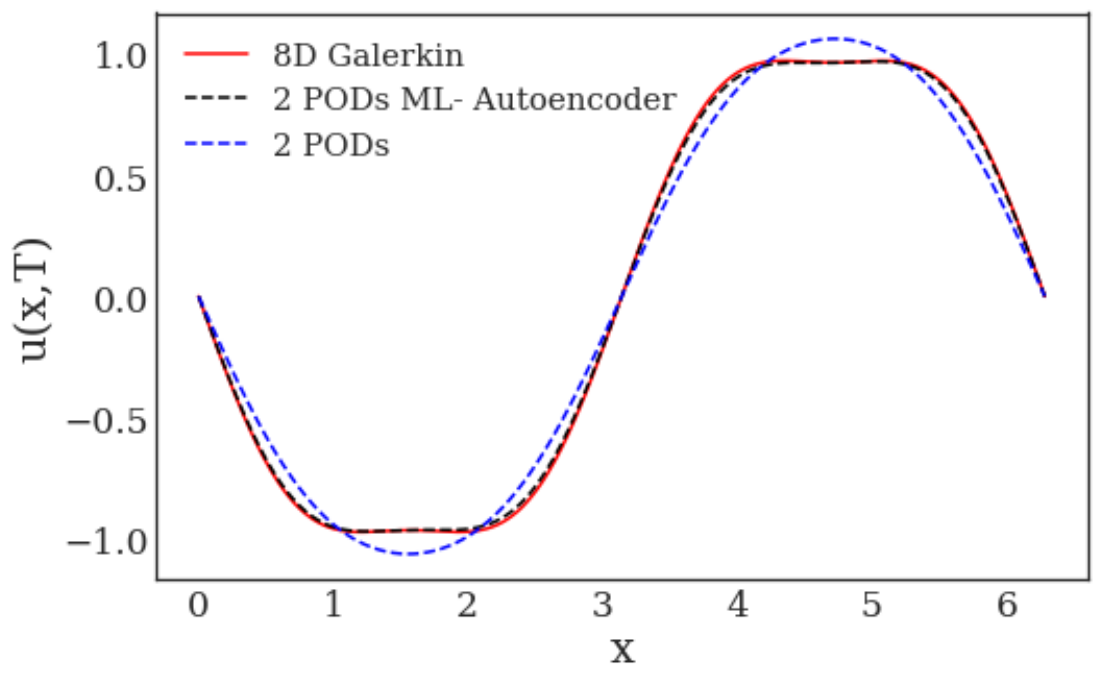}
    \caption{Reconstructed solution of the Chafee-Infante equation, at time $T=5$. The red line represents the ground truth, 8D Galerkin solution, which corresponds to 3 POD coefficients. The black broken line corresponds to the data-driven post-processed solution of the evolution of 2 POD coefficients. The uncorrected solution derived by 2 POD coefficients, is also depicted with a blue broken line.}
    \label{fig:POD2}
\end{figure}

The time-evolution law of the ODE for the two leading POD coefficients, is then learned from data. This is achieved using a feed forward neural network consisting of two hidden layers with 20 neurons each. The $tanh$ activation function is implemented and the mean squared error is used as a loss function. The learned ODE is integrated with a Runge-Kutta solver over time $T=5$. From the values of the two POD coefficients at the final time-step, the latent variables are then inferred using an appropriately trained neural network. Then, the decoder of the autoencoder is used to recover the entire set of POD coefficients. It is then possible to ``lift" from POD space to the sine coefficients and reconstruct the solution: the solution reconstructed using 3 ML-derived terms compares very favorably to the ground truth solution in Fig.\ref{fig:POD2}. For reference, in the same figure, the solution reconstructed from only the 2 POD coefficients, is also included.

\begin{figure}[ht!]
\centering
    \begin{subfigure}[b]{0.4\textwidth}
        \centering
        \includegraphics[width=\textwidth]{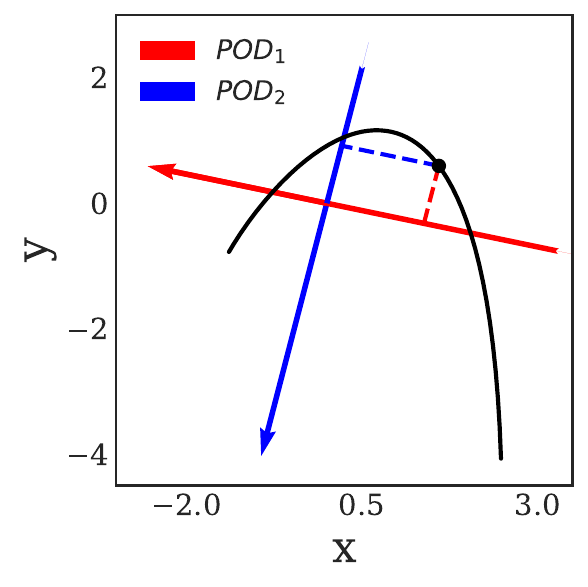}
        \caption{}
    \end{subfigure}
    \begin{subfigure}[b]{0.4\textwidth}
        \centering        \includegraphics[width=\textwidth]{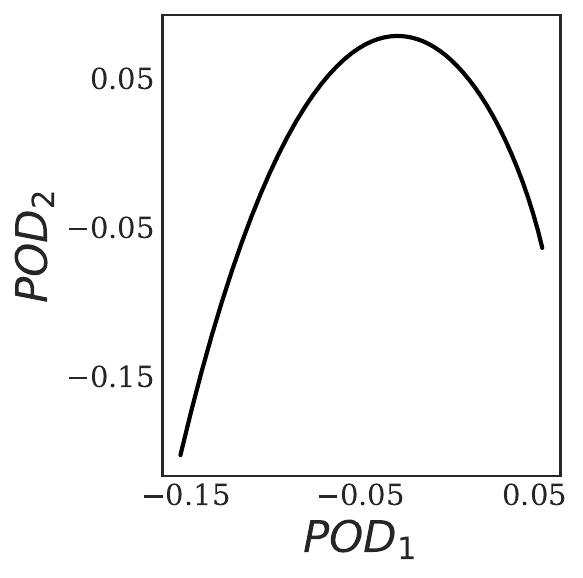}   
        \caption{}
    \end{subfigure}
    \caption{(a) A data set sampled from the singularly perturbed system of ODEs is shown with a black solid line. The span of the first POD mode ($POD_1$) is shown with a red vector and the span of the second POD mode ($POD_2$) is shown with a blue vector. The projection of a data point (black solid circle) to $POD_1$ and $POD_2$ is depicted. (b) The components of the first POD vector ($POD_1$) versus the components of the second POD vector ($POD_2$). $POD_2$ can be seen as a quadratic function of $POD_1$.}
    \label{fig:POD1_POD2_example}
\end{figure}

\section{Conclusions}
In conclusion, this study has attempted to bridge theoretical approaches to reduced order modeling of dynamical systems (theoretically, closed form,  approximations of AIMs and AIFs) with appropriately derived data-driven workflows. The data in question may consist of either (a) theoretical parametrizations of the IM (here sine coefficients) or (b) equally possibly, data-driven parametrizations (POD coefficients, autoencoder latent variables, manifold learning Diffusion Map coordinates). The use of machine learning techniques, specifically autoencoders and Diffusion Maps, allows for accurate and efficient modeling of high-dimensional systems while overcoming the limitations of traditional post-processing Galerkin methods.

Moreover, the proposed approach has demonstrated promising results in scenarios where the low-dimensional ROM significantly deviates from the correct long-term dynamics, which was previously challenging to address with post-processing Galerkin techniques. The introduction of a ``Gray-Box" model that adds a correction to the truncated Galerkin helps it regain its accuracy; it then allows for post-processing steps to recover even higher levels of accuracy, in ambient space.

Overall, this work contributes to the growing body of literature on data-driven reduced order modeling techniques for dynamical systems and provides a valuable alternative to traditional post-processing Galerkin methods. The proposed workflows have the potential to significantly improve the accuracy and efficiency of reduced order models, which has important implications for a wide range of applications, including but not limited to, aerospace engineering, biomedical engineering, and climate modeling.

A promising future direction of our current work for the construction of reduced-order models is the combination of data-driven techniques with physics-based techniques. The work of R. Geelen et al. \cite{geelen2023operator}, in which the parameterization of the data is achieved by combining linear subspaces - spanned by the first few POD vectors - and quadratic components, is the most pertinent to this direction. One could express the dynamics in terms of the first few POD vectors and use the quadratic correction only as a post-processing step to obtain a more accurate reconstruction at the end of the integration. The ability to find a quadratic correction could provide improved explainability to the post-processing step, that we lose by learning a black-box post-processing step in our current work. A visualizable example is shown in Figure \ref{fig:POD1_POD2_example} where the 2-dimensional singularly perturbed system ($\dot{x} = 2-x-y; \dot{y} = 1/{\epsilon}(x-y)$) was used to sample data. For this example one could \textit{write} the dynamics in terms of $POD_1$ and express the correction from $POD_2=f(POD_1)$ through a quadratic correction since $POD_2$ can be seen as a quadratic function of $POD_1$ (Figure \ref{fig:POD1_POD2_example}(b)). 

\subsection{Acknowledgements}
The motivation for this work comes in part from initial efforts on reduced modeling of multiphase flows, as part of CML's Thesis. I.G.K. acknowledges partial support from the US AFOSR FA9550-21-0317 and the US Department of Energy SA22-0052-S001.
C.M.L. received the support of a ``la Caixa'' Foundation Fellowship (ID
100010434), code LCF/BQ/AA19/11720048.
The research of EST was made possible by NPRP grant \#S-0207-200290 from 
the Qatar National Research Fund (a member of Qatar Foundation), and is 
based upon work supported by King Abdullah University of Science and 
Technology (KAUST) Office of Sponsored Research (OSR) under Award No. 
OSR-2020-CRG9-4336. The work of EST has also benefited from the 
inspiring environment of the CRC 1114 ``Scaling Cascades in Complex 
Systems'', Project Number 235221301, Project A02, funded by Deutsche 
Forschungsgemeinschaft (DFG). For the purpose of open access, EST has applied a Creative Commons Attribution (CC BY) licence to any Author Accepted Manuscript version arising from this submission.

\bibliographystyle{plainnat}
\bibliography{arxiv_template}

\begin{thebibliography}{63}
\providecommand{\natexlab}[1]{#1}
\providecommand{\url}[1]{\texttt{#1}}
\expandafter\ifx\csname urlstyle\endcsname\relax
  \providecommand{\doi}[1]{doi: #1}\else
  \providecommand{\doi}{doi: \begingroup \urlstyle{rm}\Url}\fi

\bibitem[Abadi et~al.(2015)Abadi, Agarwal, Barham, Brevdo, Chen, Citro, Corrado, Davis, Dean, Devin, Ghemawat, Goodfellow, Harp, Irving, Isard, Jia, Jozefowicz, Kaiser, Kudlur, Levenberg, Man\'{e}, Monga, Moore, Murray, Olah, Schuster, Shlens, Steiner, Sutskever, Talwar, Tucker, Vanhoucke, Vasudevan, Vi\'{e}gas, Vinyals, Warden, Wattenberg, Wicke, Yu, and Zheng]{tensorflow2015-whitepaper}
Mart\'{i}n Abadi, Ashish Agarwal, Paul Barham, Eugene Brevdo, Zhifeng Chen, Craig Citro, Greg~S. Corrado, Andy Davis, Jeffrey Dean, Matthieu Devin, Sanjay Ghemawat, Ian Goodfellow, Andrew Harp, Geoffrey Irving, Michael Isard, Yangqing Jia, Rafal Jozefowicz, Lukasz Kaiser, Manjunath Kudlur, Josh Levenberg, Dandelion Man\'{e}, Rajat Monga, Sherry Moore, Derek Murray, Chris Olah, Mike Schuster, Jonathon Shlens, Benoit Steiner, Ilya Sutskever, Kunal Talwar, Paul Tucker, Vincent Vanhoucke, Vijay Vasudevan, Fernanda Vi\'{e}gas, Oriol Vinyals, Pete Warden, Martin Wattenberg, Martin Wicke, Yuan Yu, and Xiaoqiang Zheng.
\newblock {TensorFlow}: Large-scale machine learning on heterogeneous systems, 2015.
\newblock URL \url{https://www.tensorflow.org/}.
\newblock Software available from tensorflow.org.

\bibitem[Adrover et~al.(2002)Adrover, Continillo, Crescitelli, Giona, and Russo]{adrover2002construction}
A~Adrover, G~Continillo, S~Crescitelli, M~Giona, and L~Russo.
\newblock Construction of approximate inertial manifold by decimation of collocation equations of distributed parameter systems.
\newblock \emph{Computers \& chemical engineering}, 26\penalty0 (1):\penalty0 113--123, 2002.

\bibitem[Akram et~al.(2020)Akram, Hassanaly, and Raman]{akram2020priori}
Maryam Akram, Malik Hassanaly, and Venkat Raman.
\newblock A priori analysis of reduced description of dynamical systems using approximate inertial manifolds.
\newblock \emph{Journal of Computational Physics}, 409:\penalty0 109344, 2020.

\bibitem[Alekseenko et~al.(1985)Alekseenko, Nakoryakov, and Pokusaev]{alekseenko1985wave}
SV~Alekseenko, VE~Nakoryakov, and BG~Pokusaev.
\newblock Wave formation on vertical falling liquid films.
\newblock \emph{International Journal of Multiphase Flow}, 11\penalty0 (5):\penalty0 607--627, 1985.

\bibitem[Anirudh et~al.(2020)Anirudh, Thiagarajan, Bremer, and Spears]{anirudh2020improved}
Rushil Anirudh, Jayaraman~J Thiagarajan, Peer-Timo Bremer, and Brian~K Spears.
\newblock Improved surrogates in inertial confinement fusion with manifold and cycle consistencies.
\newblock \emph{Proceedings of the National Academy of Sciences}, 117\penalty0 (18):\penalty0 9741--9746, 2020.

\bibitem[Bar-Sinai et~al.(2019)Bar-Sinai, Hoyer, Hickey, and Brenner]{bar2019learning}
Yohai Bar-Sinai, Stephan Hoyer, Jason Hickey, and Michael~P Brenner.
\newblock Learning data-driven discretizations for partial differential equations.
\newblock \emph{Proceedings of the National Academy of Sciences}, 116\penalty0 (31):\penalty0 15344--15349, 2019.

\bibitem[Benner et~al.(2015)Benner, Gugercin, and Willcox]{benner2015survey}
Peter Benner, Serkan Gugercin, and Karen Willcox.
\newblock A survey of projection-based model reduction methods for parametric dynamical systems.
\newblock \emph{SIAM review}, 57\penalty0 (4):\penalty0 483--531, 2015.

\bibitem[Brunton et~al.(2016)Brunton, Proctor, and Kutz]{brunton2016discovering}
Steven~L Brunton, Joshua~L Proctor, and J~Nathan Kutz.
\newblock Discovering governing equations from data by sparse identification of nonlinear dynamical systems.
\newblock \emph{Proceedings of the National Academy of Sciences}, 113\penalty0 (15):\penalty0 3932--3937, 2016.

\bibitem[Chang(1986{\natexlab{a}})]{chang1986nonlinear}
Hsueh-Chia Chang.
\newblock Nonlinear waves on liquid film surfaces—i. flooding in a vertical tube.
\newblock \emph{Chemical Engineering Science}, 41\penalty0 (10):\penalty0 2463--2476, 1986{\natexlab{a}}.

\bibitem[Chang(1986{\natexlab{b}})]{chang1986traveling}
Hsueh-Chia Chang.
\newblock Traveling waves on fluid interfaces: normal form analysis of the {K}uramoto--{S}ivashinsky equation.
\newblock \emph{The Physics of Fluids}, 29\penalty0 (10):\penalty0 3142--3147, 1986{\natexlab{b}}.

\bibitem[Chorin and Lu(2015)]{chorin2015discrete}
Alexandre~J Chorin and Fei Lu.
\newblock Discrete approach to stochastic parametrization and dimension reduction in nonlinear dynamics.
\newblock \emph{Proceedings of the National Academy of Sciences}, 112\penalty0 (32):\penalty0 9804--9809, 2015.

\bibitem[Coifman and Lafon(2006{\natexlab{a}})]{coifman2006geometric}
Ronald~R Coifman and St{\'e}phane Lafon.
\newblock Geometric harmonics: a novel tool for multiscale out-of-sample extension of empirical functions.
\newblock \emph{Applied and Computational Harmonic Analysis}, 21\penalty0 (1):\penalty0 31--52, 2006{\natexlab{a}}.

\bibitem[Coifman and Lafon(2006{\natexlab{b}})]{r19}
Ronald~R Coifman and St{\'e}phane Lafon.
\newblock Diffusion maps.
\newblock \emph{Applied and computational harmonic analysis}, 21\penalty0 (1):\penalty0 5--30, 2006{\natexlab{b}}.

\bibitem[Coifman et~al.(2008)Coifman, Kevrekidis, Lafon, Maggioni, and Nadler]{r21}
Ronald~R Coifman, Ioannis~G Kevrekidis, St{\'e}phane Lafon, Mauro Maggioni, and Boaz Nadler.
\newblock Diffusion maps, reduction coordinates, and low dimensional representation of stochastic systems.
\newblock \emph{Multiscale Modeling \& Simulation}, 7\penalty0 (2):\penalty0 842--864, 2008.

\bibitem[De~Jes{\'u}s and Graham(2023)]{de2023data}
Carlos E~P{\'e}rez De~Jes{\'u}s and Michael~D Graham.
\newblock Data-driven low-dimensional dynamic model of kolmogorov flow.
\newblock \emph{Physical Review Fluids}, 8\penalty0 (4):\penalty0 044402, 2023.

\bibitem[Dsilva et~al.(2018)Dsilva, Talmon, Coifman, and Kevrekidis]{r25}
Carmeline~J Dsilva, Ronen Talmon, Ronald~R Coifman, and Ioannis~G Kevrekidis.
\newblock Parsimonious representation of nonlinear dynamical systems through manifold learning: A chemotaxis case study.
\newblock \emph{Applied and Computational Harmonic Analysis}, 44\penalty0 (3):\penalty0 759--773, 2018.

\bibitem[Evangelou et~al.(2022)Evangelou, Dietrich, Chiavazzo, Lehmberg, Meila, and Kevrekidis]{evangelou2022double}
Nikolaos Evangelou, Felix Dietrich, Eliodoro Chiavazzo, Daniel Lehmberg, Marina Meila, and Ioannis~G Kevrekidis.
\newblock Double diffusion maps and their latent harmonics for scientific computations in latent space.
\newblock \emph{arXiv preprint arXiv:2204.12536}, 2022.

\bibitem[Evangelou et~al.(2023)Evangelou, Dietrich, Chiavazzo, Lehmberg, Meila, and Kevrekidis]{evangelou2023double}
Nikolaos Evangelou, Felix Dietrich, Eliodoro Chiavazzo, Daniel Lehmberg, Marina Meila, and Ioannis~G Kevrekidis.
\newblock Double diffusion maps and their latent harmonics for scientific computations in latent space.
\newblock \emph{Journal of Computational Physics}, 485:\penalty0 112072, 2023.

\bibitem[Foias et~al.(1988{\natexlab{a}})Foias, Jolly, Kevrekidis, Sell, and Titi]{foias1988computation}
C~Foias, MS~Jolly, IG~Kevrekidis, George~R Sell, and ES~Titi.
\newblock On the computation of inertial manifolds.
\newblock \emph{Physics Letters A}, 131\penalty0 (7-8):\penalty0 433--436, 1988{\natexlab{a}}.

\bibitem[Foias et~al.(1988{\natexlab{b}})Foias, Sell, and Temam]{foias1988inertial}
Ciprian Foias, George~R Sell, and Roger Temam.
\newblock Inertial manifolds for nonlinear evolutionary equations.
\newblock \emph{Journal of {d}ifferential {e}quations}, 73\penalty0 (2):\penalty0 309--353, 1988{\natexlab{b}}.

\bibitem[Foias et~al.(1989)Foias, Sell, and Titi]{foias1989exponential}
Ciprian Foias, George~R Sell, and Edriss~S Titi.
\newblock Exponential tracking and approximation of inertial manifolds for dissipative nonlinear equations.
\newblock \emph{Journal of Dynamics and Differential Equations}, 1:\penalty0 199--244, 1989.

\bibitem[Garc{\'\i}a-Archilla and Titi(1999)]{garcia1999postprocessing}
Bosco Garc{\'\i}a-Archilla and Edriss~S Titi.
\newblock Postprocessing the {G}alerkin method: the finite-element case.
\newblock \emph{SIAM Journal on Numerical Analysis}, 37\penalty0 (2):\penalty0 470--499, 1999.

\bibitem[Garc{\'\i}a-Archilla et~al.(1998)Garc{\'\i}a-Archilla, Novo, and Titi]{garcia1998postprocessing}
Bosco Garc{\'\i}a-Archilla, Julia Novo, and Edriss~S Titi.
\newblock Postprocessing the {G}alerkin method: a novel approach to approximate inertial manifolds.
\newblock \emph{SIAM Journal on Numerical Analysis}, 35\penalty0 (3):\penalty0 941--972, 1998.

\bibitem[Garc{\'\i}a-Archilla et~al.(1999)Garc{\'\i}a-Archilla, Novo, and Titi]{garcia1999approximate}
Bosco Garc{\'\i}a-Archilla, Julia Novo, and Edriss Titi.
\newblock An approximate inertial manifolds approach to postprocessing the galerkin method for the navier-stokes equations.
\newblock \emph{Mathematics of Computation}, 68\penalty0 (227):\penalty0 893--911, 1999.

\bibitem[Gear et~al.(2011)Gear, Kevrekidis, and Sonday]{gear2011slow}
Charles~William Gear, IG~Kevrekidis, and BE~Sonday.
\newblock Slow manifold integration on a diffusion map parameterization.
\newblock In \emph{AIP Conference Proceedings}, volume 1389, pages 13--16. American Institute of Physics, 2011.

\bibitem[Geelen et~al.(2023)Geelen, Wright, and Willcox]{geelen2023operator}
Rudy Geelen, Stephen Wright, and Karen Willcox.
\newblock Operator inference for non-intrusive model reduction with quadratic manifolds.
\newblock \emph{Computer Methods in Applied Mechanics and Engineering}, 403:\penalty0 115717, 2023.

\bibitem[Guermond and Prudhomme(2008)]{guermond2008fully}
J-L Guermond and Serge Prudhomme.
\newblock A fully discrete nonlinear {G}alerkin method for the {3D} {N}avier--{S}tokes equations.
\newblock \emph{Numerical Methods for Partial Differential Equations: An International Journal}, 24\penalty0 (3):\penalty0 759--775, 2008.

\bibitem[Jauberteau et~al.(1990)Jauberteau, Rosier, and Temam]{jauberteau1990nonlinear}
F~Jauberteau, C~Rosier, and R~Temam.
\newblock A nonlinear {G}alerkin method for the navier-stokes equations.
\newblock \emph{Computer Methods in Applied Mechanics and Engineering}, 80\penalty0 (1-3):\penalty0 245--260, 1990.

\bibitem[Jolly(1989)]{jolly1989explicit}
Michael~S Jolly.
\newblock Explicit construction of an inertial manifold for a reaction diffusion equation.
\newblock \emph{Journal of Differential Equations}, 78\penalty0 (2):\penalty0 220--261, 1989.

\bibitem[Jolly et~al.(1990)Jolly, Kevrekidis, and Titi]{jolly1990approximate}
Michael~S Jolly, IG~Kevrekidis, and Edriss~S Titi.
\newblock Approximate inertial manifolds for the {K}uramoto-{S}ivashinsky equation: analysis and computations.
\newblock \emph{Physica D: Nonlinear Phenomena}, 44\penalty0 (1-2):\penalty0 38--60, 1990.

\bibitem[Jolly et~al.(1991)Jolly, Kevrekidis, and Titi]{jolly1991preserving}
MS~Jolly, IG~Kevrekidis, and ES~Titi.
\newblock Preserving dissipation in approximate inertial forms for the {K}uramoto-{S}ivashinsky equation.
\newblock \emph{Journal of Dynamics and Differential Equations}, 3:\penalty0 179--197, 1991.

\bibitem[Kang et~al.(2015)Kang, Zhang, Ren, and Lei]{kang2015nonlinear}
Wei Kang, Jia-Zhong Zhang, Sheng Ren, and Peng-Fei Lei.
\newblock Nonlinear {G}alerkin method for low-dimensional modeling of fluid dynamic system using pod modes.
\newblock \emph{Communications in Nonlinear Science and Numerical Simulation}, 22\penalty0 (1-3):\penalty0 943--952, 2015.

\bibitem[Kevrekidis et~al.(1990)Kevrekidis, Nicolaenko, and Scovel]{Kevrekidis1990}
Ioannis~G. Kevrekidis, Basil Nicolaenko, and James~C. Scovel.
\newblock Back in the saddle again: A computer assisted study of the {K}uramoto-{S}ivashinsky equation.
\newblock \emph{SIAM Journal on Applied Mathematics}, 50\penalty0 (3):\penalty0 760--790, 1990.
\newblock ISSN 00361399.
\newblock URL \url{http://www.jstor.org/stable/2101886}.

\bibitem[Kramer(1991)]{kramer1991nonlinear}
Mark~A Kramer.
\newblock Nonlinear principal component analysis using autoassociative neural networks.
\newblock \emph{AIChE journal}, 37\penalty0 (2):\penalty0 233--243, 1991.

\bibitem[Krischer et~al.(1993)Krischer, Rico-Martínez, Kevrekidis, Rotermund, Ertl, and Hudson]{KrischerKevrekidisErtl1993}
K.~Krischer, R.~Rico-Martínez, I.~G. Kevrekidis, H.~H. Rotermund, G.~Ertl, and J.~L. Hudson.
\newblock Model identification of a spatiotemporally varying catalytic reaction.
\newblock \emph{AIChE Journal}, 39\penalty0 (1):\penalty0 89--98, 1993.
\newblock \doi{https://doi.org/10.1002/aic.690390110}.
\newblock URL \url{https://aiche.onlinelibrary.wiley.com/doi/abs/10.1002/aic.690390110}.

\bibitem[Kuramoto and Tsuzuki(1976)]{kuramoto1976persistent}
Yoshiki Kuramoto and Toshio Tsuzuki.
\newblock Persistent propagation of concentration waves in dissipative media far from thermal equilibrium.
\newblock \emph{Progress of theoretical physics}, 55\penalty0 (2):\penalty0 356--369, 1976.

\bibitem[Lee and Carlberg(2020)]{lee2020model}
Kookjin Lee and Kevin~T Carlberg.
\newblock Model reduction of dynamical systems on nonlinear manifolds using deep convolutional autoencoders.
\newblock \emph{Journal of Computational Physics}, 404:\penalty0 108973, 2020.

\bibitem[Linot and Graham(2020)]{linot2020deep}
Alec~J Linot and Michael~D Graham.
\newblock Deep learning to discover and predict dynamics on an inertial manifold.
\newblock \emph{Physical Review E}, 101\penalty0 (6):\penalty0 062209, 2020.

\bibitem[Linot and Graham(2022)]{linot2022data}
Alec~J Linot and Michael~D Graham.
\newblock Data-driven reduced-order modeling of spatiotemporal chaos with neural ordinary differential equations.
\newblock \emph{Chaos: An Interdisciplinary Journal of Nonlinear Science}, 32\penalty0 (7):\penalty0 073110, 2022.

\bibitem[Linot et~al.(2023)Linot, Zeng, and Graham]{linot2023turbulence}
Alec~J Linot, Kevin Zeng, and Michael~D Graham.
\newblock Turbulence control in plane couette flow using low-dimensional neural ode-based models and deep reinforcement learning.
\newblock \emph{International Journal of Heat and Fluid Flow}, 101:\penalty0 109139, 2023.

\bibitem[Lu et~al.(2017)Lu, Lin, and Chorin]{lu2017data}
Fei Lu, Kevin~K Lin, and Alexandre~J Chorin.
\newblock Data-based stochastic model reduction for the {K}uramoto--{S}ivashinsky equation.
\newblock \emph{Physica D: Nonlinear Phenomena}, 340:\penalty0 46--57, 2017.

\bibitem[Margolin et~al.(2003)Margolin, Titi, and Wynne]{margolin2003postprocessing}
Len~G Margolin, Edriss~S Titi, and Shannon Wynne.
\newblock The postprocessing {G}alerkin and nonlinear {G}alerkin methods---a truncation analysis point of view.
\newblock \emph{SIAM Journal on Numerical Analysis}, 41\penalty0 (2):\penalty0 695--714, 2003.

\bibitem[Marion and Temam(1990)]{marion1990nonlinear}
Martine Marion and R~Temam.
\newblock Nonlinear {G}alerkin methods: the finite elements case.
\newblock \emph{Numerische Mathematik}, 57:\penalty0 205--226, 1990.

\bibitem[Marsden et~al.(1993)Marsden, Hoffman, et~al.]{marsden1993elementary}
Jerrold~E Marsden, Michael~J Hoffman, et~al.
\newblock \emph{Elementary classical analysis}.
\newblock Macmillan, 1993.

\bibitem[Martin-Linares et~al.(2023)Martin-Linares, Psarellis, Karapetsas, Koronaki, and Kevrekidis]{martin2023physics}
Cristina~P Martin-Linares, Yorgos~M Psarellis, Georgios Karapetsas, Eleni~D Koronaki, and Ioannis~G Kevrekidis.
\newblock Physics-agnostic and physics-infused machine learning for thin films flows: modeling, and predictions from small data.
\newblock \emph{arXiv preprint arXiv:2301.12508}, 2023.

\bibitem[McQuarrie et~al.(2021)McQuarrie, Huang, and Willcox]{mcquarrie2021data}
Shane~A McQuarrie, Cheng Huang, and Karen~E Willcox.
\newblock Data-driven reduced-order models via regularised operator inference for a single-injector combustion process.
\newblock \emph{Journal of the Royal Society of New Zealand}, 51\penalty0 (2):\penalty0 194--211, 2021.

\bibitem[Nadler et~al.(2006)Nadler, Lafon, Coifman, and Kevrekidis]{r20}
Boaz Nadler, St{\'e}phane Lafon, Ronald~R Coifman, and Ioannis~G Kevrekidis.
\newblock Diffusion maps, spectral clustering and reaction coordinates of dynamical systems.
\newblock \emph{Applied and Computational Harmonic Analysis}, 21\penalty0 (1):\penalty0 113--127, 2006.

\bibitem[Qian et~al.(2022)Qian, Farcas, and Willcox]{qian2022reduced}
Elizabeth Qian, Ionut-Gabriel Farcas, and Karen Willcox.
\newblock Reduced operator inference for nonlinear partial differential equations.
\newblock \emph{SIAM Journal on Scientific Computing}, 44\penalty0 (4):\penalty0 A1934--A1959, 2022.

\bibitem[Raissi et~al.(2019)Raissi, Perdikaris, and Karniadakis]{raissi2019physics}
Maziar Raissi, Paris Perdikaris, and George~E Karniadakis.
\newblock Physics-informed neural networks: A deep learning framework for solving forward and inverse problems involving nonlinear partial differential equations.
\newblock \emph{Journal of Computational physics}, 378:\penalty0 686--707, 2019.

\bibitem[Rico-Martinez et~al.(1992)Rico-Martinez, Krischer, Kevrekidis, Kube, and Hudson]{rico1992discrete}
Ramiro Rico-Martinez, K~Krischer, IG~Kevrekidis, MC~Kube, and JL~Hudson.
\newblock Discrete-vs. continuous-time nonlinear signal processing of cu electrodissolution data.
\newblock \emph{Chemical Engineering Communications}, 118\penalty0 (1):\penalty0 25--48, 1992.

\bibitem[Shen(1990)]{shen1990long}
Jie Shen.
\newblock Long time stability and convergence for fully discrete nonlinear {G}alerkin methods.
\newblock \emph{Applicable Analysis}, 38\penalty0 (4):\penalty0 201--229, 1990.

\bibitem[Shvartsman and Kevrekidis(1998)]{shvartsman1998nonlinear}
Stanislav~Y Shvartsman and Ioannis~G Kevrekidis.
\newblock Nonlinear model reduction for control of distributed systems: A computer-assisted study.
\newblock \emph{AIChE Journal}, 44\penalty0 (7):\penalty0 1579--1595, 1998.

\bibitem[Sivashinsky(1977)]{sivashinsky1977nonlinear}
Gregory~I Sivashinsky.
\newblock Nonlinear analysis of hydrodynamic instability in laminar flames—i. derivation of basic equations.
\newblock \emph{Acta astronautica}, 4\penalty0 (11):\penalty0 1177--1206, 1977.

\bibitem[Sonday(2011)]{sonday2011systematic}
Benjamin Sonday.
\newblock \emph{Systematic model reduction for complex systems through data mining and dimensionality reduction}.
\newblock Princeton University, 2011.

\bibitem[Sonday et~al.(2010)Sonday, Singer, Gear, and Kevrekidis]{sonday2010manifold}
Benjamin~E Sonday, Amit Singer, C~William Gear, and Ioannis~G Kevrekidis.
\newblock Manifold learning techniques and model reduction applied to dissipative pdes.
\newblock \emph{arXiv preprint arXiv:1011.5197}, 2010.

\bibitem[Temam(1989{\natexlab{a}})]{temam1989inertial}
R~Temam.
\newblock Do inertial manifolds apply to turbulence?
\newblock \emph{Physica D: Nonlinear Phenomena}, 37\penalty0 (1-3):\penalty0 146--152, 1989{\natexlab{a}}.

\bibitem[Temam(1989{\natexlab{b}})]{temam1989induced}
Roger Temam.
\newblock Induced trajectories and approximate inertial manifolds.
\newblock \emph{ESAIM: Mathematical Modelling and Numerical Analysis}, 23\penalty0 (3):\penalty0 541--561, 1989{\natexlab{b}}.

\bibitem[Theodoropoulos et~al.(2000)Theodoropoulos, Kevrekidis, and Mountziaris]{shvartsman2000order}
C~Theodoropoulos, IG~Kevrekidis, and TJ~Mountziaris.
\newblock Order reduction for nonlinear dynamic models of distributed reacting systems.
\newblock \emph{Journal of Process Control}, 10\penalty0 (2-3):\penalty0 177--184, 2000.

\bibitem[Titi(1990)]{titi1990approximate}
Edriss~S Titi.
\newblock On approximate inertial manifolds to the navier-stokes equations.
\newblock \emph{Journal of Mathematical Analysis and Applications}, 149\penalty0 (2):\penalty0 540--557, 1990.

\bibitem[Wahlbin(2006)]{wahlbin2006superconvergence}
Lars Wahlbin.
\newblock \emph{Superconvergence in Galerkin finite element methods}.
\newblock Springer, 2006.

\bibitem[Zastrow et~al.(2023)Zastrow, Chaudhuri, Willcox, Ashley, and Henson]{zastrow2023data}
Benjamin~G Zastrow, Anirban Chaudhuri, Karen~E Willcox, Anthony~S Ashley, and Michael~C Henson.
\newblock Data-driven model reduction via operator inference for coupled aeroelastic flutter.
\newblock In \emph{AIAA SCITECH 2023 Forum}, page 0330, 2023.

\bibitem[Zeng and Graham(2023)]{zeng2023autoencoders}
Kevin Zeng and Michael~D. Graham.
\newblock Autoencoders for discovering manifold dimension and coordinates in data from complex dynamical systems, 2023.

\bibitem[Zeng et~al.(2022)Zeng, Linot, and Graham]{zeng2022data}
Kevin Zeng, Alec~J Linot, and Michael~D Graham.
\newblock Data-driven control of spatiotemporal chaos with reduced-order neural ode-based models and reinforcement learning.
\newblock \emph{Proceedings of the Royal Society A}, 478\penalty0 (2267):\penalty0 20220297, 2022.

\end{thebibliography}

\appendix
\section{Appendix}
\subsubsection{Post-Processing Galerkin for the finite element Method}
\label{sec:PPE_FE}

Let $X$ be the phase space of a nonlinear dissipative evolution equation of the form
\begin{equation*}
    \frac{du}{dt}+\nu Au=F(u).
\end{equation*}
\noindent    
Let $X_H$ be a finite dimensional (e.g. finite element) space of spatial scale $H$ with $\mathcal{P}_H$: $X\rightarrow X_H$ an orthogonal projection. The Galerkin approximate solution $u_H \in X_H$ satisfies the equation
\begin{equation*}
\frac{du_H}{dt} + \nu \mathcal{P}_H A u_H=\mathcal{P}_H F(u_H), \;\; \text{for} \;  t\in [0,T].
\end{equation*}

\noindent
Therefore, for a given Galerkin solution $u_H$ and its time-derivative $\frac{du_H}{dt}$ over the interval $[0,T]$, the post-processing Galekin solution is a function $v \in X$ (notice it is not in the complement of $X_H$, i.e. not in $X \ominus X_H$, but in $X$, so $v$ involves both coarse as well as fine spatial scales), such that $v$ satisfies

\begin{equation*}
    \nu A v = - \frac{du_H}{dt} + F(u_H) \mid_{t=T}.
\end{equation*}

\noindent
The right-hand side is a given function at time $t=T$, and $v$ solves a linear elliptic equation. However, in practice we solve an approximation of $v$ say $\tilde{v} \in X_h$, where $h \ll H$, and $X_h$ is a finer finite element space
\begin{equation*}
    \nu \mathcal{P}_h A \tilde{v}=\mathcal{P}_h \bigg(-\frac{du_H}{dt}+f(u_H) \bigg)\mid_{t=T}.
\end{equation*}

\subsubsection{Euler-Galerkin algorithm applied to Chafee-Infante PDE}
\label{sec:Euler_Garlerkin_Chafee_Infante}
The implementation of the Euler-Galerkin algorithm described in Sec. \ref{sec:Euler_Galerkin} is shown here for the Chafee-Infante reaction-diffusion equation. For this PDE, as discussed in Sec. \ref{sec:Chafee_Infante}, a two-dimensional inertial manifold exists ($n=2$) parameterized by the first two sine Fourier modes $\alpha_1, \alpha_2$. By using the Galerkin projection $u(x,t) \approx \sum_{i=1}^{m=3} a_i(t)sin(ix)$ a system of three coupled ordinary differential equations is derived. The derived system of equations reads
\begin{align}
    &\dot{a_1} = -a_1\nu  + a_1 - \frac{3}{4}a^3_1  -\frac{3}{2}a_1a^2_2 - \frac{3}{4}a^2_1a_3 -\frac{3}{4}a^2_2a_3 -\frac{3}{4}a_1a^2_3 \\
    &\dot{a_2} = -4a_2\nu + a_2 - \frac{3}{2}a^2_1a_2 - \frac{3}{4}a^3_2 - \frac{3}{2}a_1a_2a_3 -\frac{3}{2}a_2a^2_3 \\
    &\dot{a_3} = -9a_3\nu + a_3 + \frac{a^3_1}{4} - \frac{3}{2}a^3_3 - \frac{3}{2}a^2_2a_3 - \frac{3}{2}a^2_1a_3 - \frac{3}{4}a_1a^2_2.
\label{eq:Odes_Scheme_a3}
\end{align}

The term $-9\alpha_3 \nu$ of the right-hand side of Equation \eqref{eq:Odes_Scheme_a3} corresponds to the diffusion term and all the other terms of the right-hand side to the reaction terms.
We take an implicit Euler step of Equation \eqref{eq:Odes_Scheme_a3} of length $\tau$ by using as initial condition $\alpha_3(t=0) = 0$. This gives us the expression 
\begin{align}
\begin{split}
    \alpha_3(\tau) = \alpha_3(0) + \tau 
 \dot{a}_3 
 \end{split}
\end{align}

By moving the diffusive term to the left-hand side and solving in terms of $\alpha_3$ we get the expression

\begin{equation}
    \alpha_3 = \frac{\tau}{(1+9\tau \nu)}\bigg(a_3 + \frac{a^3_1}{4} - \frac{3}{2}a^3_3 - \frac{3}{2}a^2_2a_3 - \frac{3}{2}a^2_1a_3 - \frac{3}{4}a_1a^2_2 \bigg).
\end{equation}

We then perform one fixed point iteration by considering $a_3 = 0$ and $\tau = 1$. This leads to the Euler-Galerkin approximation
\begin{equation}
    \label{eq:Euler_Galerkin_AIM}
    \alpha_3  = \frac{1}{4(1+9 \nu)}\bigg(\alpha_1^3 - \alpha_1 \alpha_2^2 \bigg).
\end{equation}

In our case, the Euler-Galerkin approximation in Equation \eqref{eq:Euler_Galerkin_AIM} was used as one of the post-processing schemes to correct the solution of $\hat{u}(x,T)$ computed from the truncated dynamics.

\subsection{Diffusion Maps}
\label{sec:Diffusion_Maps_algorithm}
The Diffusion Maps algorithm reveals the \textit{intrinsic geometry} of a data set $\mathbf{X} = \{\vect{x}_i\}_{i=1}^N$,  where each data point $\vect{x}_i \in \mathbb{R}^m$, by constructing a random walk on $\mathbf{X}$. This random walk is constructed by means of an affinity matrix $\mathbf{A}$ with entries computed in terms of a kernel

\begin{equation}
\label{eq:affinity_matrix}
    A_{ij} = \text{exp} \Bigg(\frac{-|| \vect{x}_i - \vect{x}_j ||^2_2}{2 \varepsilon} \Bigg)
\end{equation}

where $\varepsilon$ is a positive hyperparameter that specifies the rate of decay of the kernel (kernel bandwidth). The Gaussian kernel, Equation \eqref{eq:affinity_matrix}, is typically chosen for the construction of the affinity matrix $\mathbf{A}$.

To obtain a random walk (parametrization) of $\mathbf{X}$ regardless of the sampling density 
the normalization 
\begin{equation}
    \mathbf{K}  = \mathbf{P}^{-\alpha} \mathbf{K} \mathbf{P}^{-\alpha}, \text{ where } P_{ii} = \sum_{j=1}^N A_{ij}
\end{equation}
is applied, with $\alpha=1$ to factor out the density effects.

A second normalization,

\begin{equation}
    \tilde{\mathbf{K}} = \mathbf{D}^{-1}\mathbf{K} \text{ where } D_{ii} = \sum_{j=1}^N K_{ij},
\end{equation}

is applied to recover the row-stochastic matrix $\tilde{\mathbf{K}}$. The eigendecomposition 
of $\tilde{\mathbf{K}}$,
\begin{equation}
    \tilde{\mathbf{K}}\vect{\phi}_i = \lambda_{i} \vect{\phi}_i
\end{equation}

gives a set of eigenvectors $\vect{\phi}$ and eigenvalues $\lambda$. Proper selection of the eigenvectors that parameterize independent directions, (known as non-harmonics) is needed. This selection in practice can be achieved by using the local linear regression algorithm proposed in \cite{r25}. If the number of non-harmonic eigenvectors is smaller than the original dimension $\tilde{n} < n$ then those eigenvectors $\mathbf{\Phi} = \{\phi_1, \dots, \phi_{\tilde{n}} \}$ can provide a more parsimonious representation of the original data and thus to obtain dimensionality reduction. 

It is therefore important to discover which eigenvectors parametrize independent directions, and do not span the same direction with different frequencies (\textit{harmonics}). 

To achieve this, the local linear regression algorithm, proposed in \cite{r25} is used, according to which, each DMAP coordinate is fitted as a function of the previous ones.
To select the DMAP coordinates that are independent, the ``goodness of fit" of this functions is used: A good fit is associated with a  $\phi_k$ that is a harmonic function of the previous eigenmodes, whereas a bad fit signifies that $\phi_k$ is a new independent direction on the data manifold.
\subsection{Geometric Harmonics and Double Diffusion Maps}
\label{sec:double_dmaps}
Geometric Harmonics  \cite{coifman2006geometric} is a regression scheme \textit{traditionally} applied on a data set $\mathbf{X}$ to extend a function $f$. Extending means that we are able to evaluate the function $f$ for points ``outside'' of $\mathbf{X}$, for $x_{new} \notin \mathbf{X}$.

In our previous work \cite{evangelou2022double} we introduced a special case of Geometric Harmonics, termed Double Diffusion Maps (Double DMAPs), able to regress functions directly on the reduced Diffusion Maps coordinates $\mathbf{\Phi}$. In this case, similar to the first round of Diffusion Maps, an affinity matrix is computed

\begin{equation}
    A^*_{ij} = \text{exp}\Bigg(\frac{-|| \vect{\phi}_i - \vect{\phi}_j ||_2^2}{2\varepsilon^{*}} \Bigg).
\end{equation}

An eigendecomposition of the symmetric and positive semidefinite matrix $\mathbf{A}^*$ is then computed. From this eigendecomposition we obtain a set of orthonormal eigenvectors $\mathbf{\Psi} = \{ \psi_0, \dots, \psi_{N-1} \}$ ranked with their non-negative eigenvalues  $\vect{\sigma} = \{\sigma_0, \dots, \sigma _{N-1} \}$. A set of those eigenvalues $S_{\delta} = \{i : \sigma_i > \delta \sigma_0$\}, where $\delta>0$, is considered as the basis in which we project and subsequently extend the function $f$. The projection of $f$ in this truncated step is given as
\begin{equation}
    f \mapsto P_{\delta}f = \sum_{i \in \Sigma_{\delta}} \langle f, \psi_i \rangle \psi_i
\end{equation}

where $\langle \cdot, \cdot \rangle$ denotes the inner product. For $\vect{\phi}_{new} \notin \mathbf{\Phi}$ we obtain $(Ef)(\vect{\phi}_{new})$ by firstly extending each eigenvector $\psi_i \in \mathbf{\Psi}$,

\begin{equation}
    \Psi_{i} (\vect{\phi}_{new}) = \sigma^{-1}_{i} \sum_{j=1}^m A(\vect{\phi}_{new}, \vect{\phi}_j)\psi_i({\phi_j}) ,
\end{equation}

where $\sigma_i$ is the $i^{\text{th}}$ eigenvalue and $\psi_i(\phi_j)$ is the $j^{\text{th}}$ component of the eigenvector $\vect{\psi}_i$. The extended eigenvectors can then used to estimate $(Ef)(\vect{\phi}_{new})$ as,

\begin{equation}
    (Ef)(\vect{\phi}_{new}) = \sum_{i \in S_{\delta}} \langle f, \psi_i \rangle \Psi_i(\phi_{new}).
 \end{equation}

\subsection{Inverse Function Theorem}
\label{sec:IFT}
Consider the vector function \hbox{$F(\vect{x}) = \vect{y}$} and assume that \hbox{$\vect{x}\in \mathbb{R}^n$}is a solution of $F$ and that $F:\mathbb{R}^n \to \mathbb{R}^n$ is differentiable . The Inverse Function Theorem \cite{marsden1993elementary} states that, if the Jacobian matrix

\begin{equation}
    \label{eqn:jacobian}
    \mathbf{J}_f(\vect{x})=\begin{pmatrix}
        \frac{\partial f_1}{\partial x_1} & \cdots & \frac{\partial f_1}{\partial x_n}\\
        \vdots & \ddots & \vdots\\
        \frac{\partial f_n}{\partial x_1} & \cdots & \frac{\partial f_n}{\partial x_n}
    \end{pmatrix}
\end{equation}

is invertible, then in a neighborhood of $\vect{x}$ and $\vect{y}$ the function $f^{-1}$ exists. This suggests a unique \textit{local} solution close to any $\vect{y}$. The Jacobian matrix is invertible if and only if its determinant is nonzero, therefore, showing that the det$(\mathbf{J}_f(\vect{x}))$ has values of a single sign guarantees that the mapping is locally invertible and thus one-to-one. 

\label{sec:SI}
\begin{figure}[ht!]
    \centering
    \includegraphics[width=\textwidth]{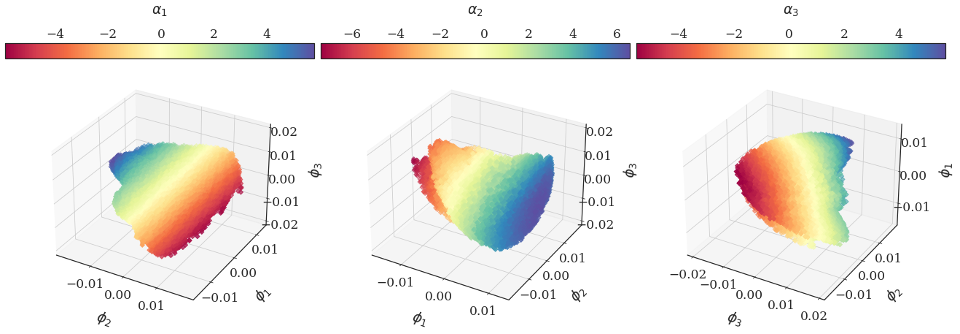}
    \caption{Diffusion Maps coordinates the parametrize the latent space: $\phi_1$, $\phi_2$ and $\phi_3$; Left: colored by sine coefficient $\alpha_1$, center: colored by sine coefficient $\alpha_2$, right: colored by sine coefficient $\alpha_3$ }
    \label{fig:ks_dmaps}
\end{figure}

\begin{figure}[ht!]
    \centering
    \includegraphics[width=0.9\textwidth]{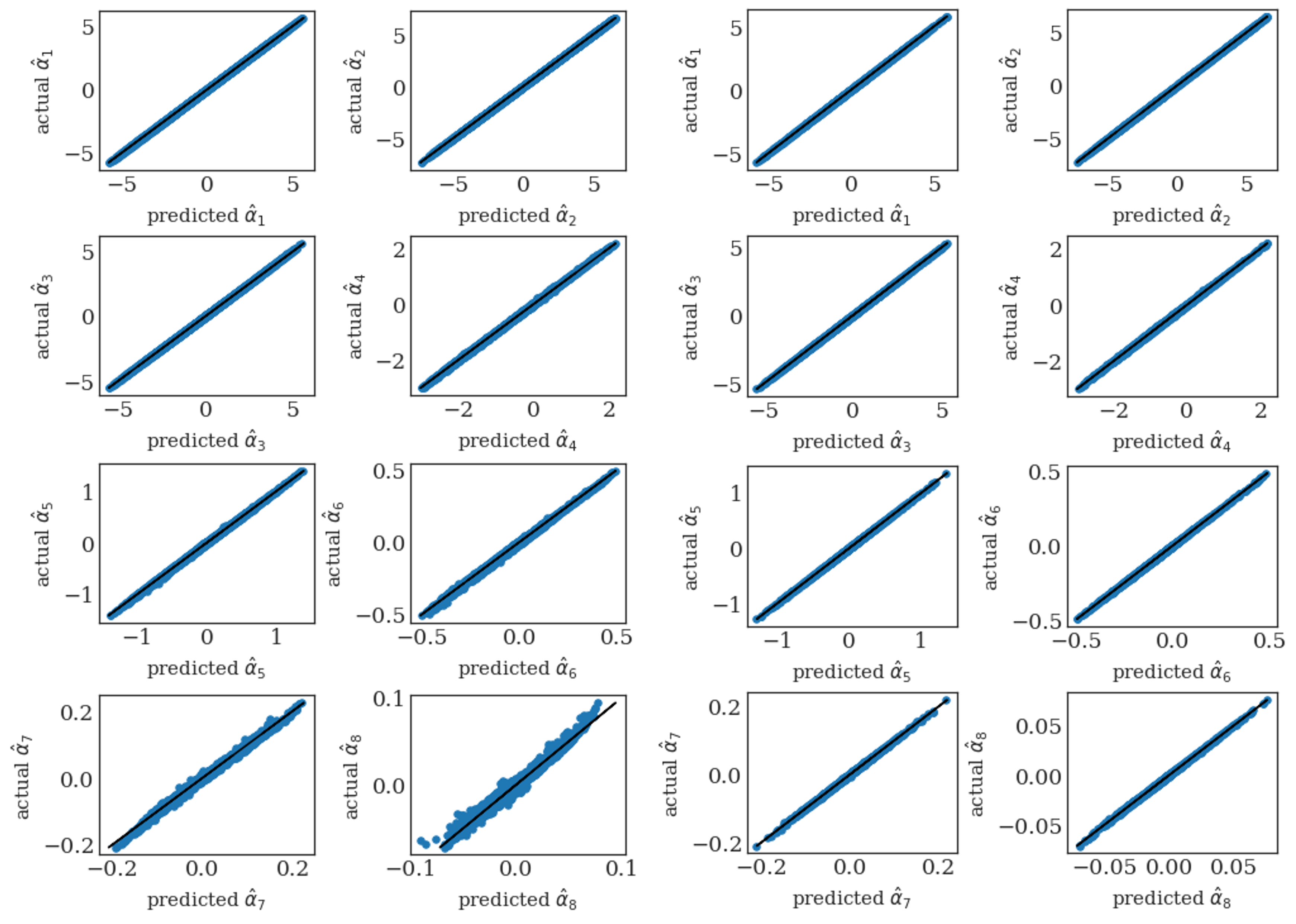}
    \caption{Reconstruction of sine coefficients; first two columns: by the autoencoder, last two columns: by Double DMAPs.}
    \label{fig:reconstructed_a}
\end{figure}

\begin{figure}[ht!]
    \centering
    \includegraphics[width=0.8\textwidth]{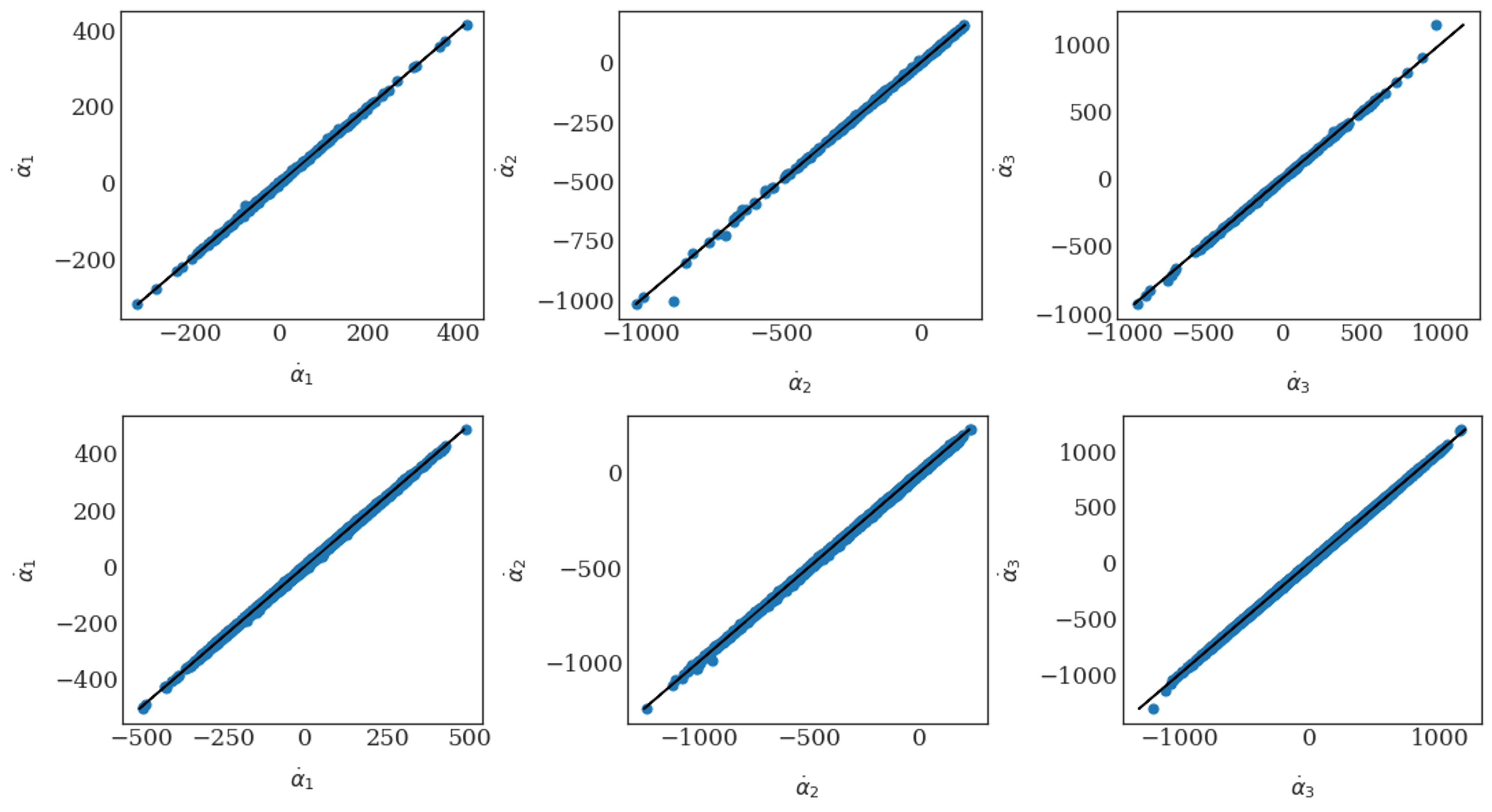}
    \caption{Performance of the neural network predicting the right-hand-side of the learned ODE of the first three coefficients. Actual versus learned $\dot{\alpha}_1$ (left), $\dot{\alpha}_2$ (center) and $\dot{\alpha}_3$ (right), from the autoencoder (top) and Double DMAPs (bottom) derived values of sine coefficients.}
    \label{fig:reconstructed_rhs}
\end{figure}

\end{document}